\documentclass[10pt]{article}
\usepackage[sumlimits]{amsmath}
\usepackage{amsmath,amsfonts,amsthm,amssymb,amscd}
\newtheorem{theo}{Theorem}[section]
\newtheorem{prop}[theo]{Proposition}
\newtheorem{lemm}[theo]{Lemma}
\newtheorem{coro}[theo]{Corollary}
\newtheorem{rem}[theo]{Remark}

\newtheorem{defi}[theo]{Definition}

\def\bbb{{\cal B}}\def\ccc{{\cal C}}
\def\fff{{\cal F}}

\def\nnn{{\cal N}} \def\ppp{{\cal P}}
  
\def\ttt{{\cal T}}

\newcommand\N{\mathbb{N}}

\newcommand\R{\mathbb{R}}

\def\1{{\bold 1}}
\def\ph{\varphi}
\def\eps{{\varepsilon}}

\title{Stochastic extensions of symbols in Wiener spaces and heat operator.}
\author{Lisette Jager}

\date{}	
		
\begin{document}

\maketitle

\abstract{The construction, in  \cite{AJN}, of a pseudodifferential calculus
 analogous to the Weyl calculus, in an infinite dimensional setting, required the introduction of convenient classes of symbols.\\ 
In this article, we proceed with the study of these classes in order to 
establish, later on, the properties that a pseudodifferential calculus 
is expected to satisfy. The introduction and the study 
of a new class are rendered necessary in view of  applications 
in QED.\\
We prove here that the symbols of both classes and the terms of their Taylor
expansions admit stochastic extensions. We define, in this infinite dimensional
 setting, a semigroup $H_t$ analogous to the heat semigroup, acting on the 
symbols belonging to both classes of symbols. The heat operator commutes 
with a second order operator similar to the Laplacian,which is its
 infinitesimal generator.
For the class defined there, we give an expansion in powers of $t$ of $H_tf$,
according to the classes of symbols.}\\


\tableofcontents

\section{Introduction}\label{intro}
This article follows \cite{AJN} where a pseudodifferential Weyl calculus 
in an infinite dimensional setting has been developped, replacing $\R^n $ by 
 a probability
space, the abstract Wiener space, denoted by $B$ (one may consult 
\cite{K} about this topic). This space is
the completion of a real, separable, infinite dimensional Hilbert space
$H$ with respect to a convenient norm (called ``measurable''), which is
different from the canonical norm of $H$.
We then have two complete spaces, one endowed with a scalar product and a
 symplectic 
form (on $H^2$), the other one, with a probability measure  $\mu_{B,t}$
 generalizing the finite dimensional Gaussian measure, the positive parameter
$t$ representing the variance. This distribution of properties compels
us to shift constantly from one space to the other, which is naturally not
 the case in the finite dimensional setting. Remark that, since the completion
 of $H$ depends on the choice of the norm, it is not unique and we shall 
take advantage of it. \\
 In \cite{AJN}, the Weyl calculus has been constructed for symbols belonging to 
a given class of symbols,  $S_m(\bbb,\eps)$, recalled in Definition 
\ref{1.3AJN}. The symbols are functions defined on the Hilbert space $H$ and
 satisfying partial differentiability conditions with respect to a fixed 
orthonormal basis $\bbb$, as well as estimates. These properties allow us to
 extend the symbols, in a certain sense, as functions defined on the Wiener
 space $B$. This is the notion of stochastic extension, recalled 
in  Definition \ref{4.4AJN} and which is generally different from
a continuity extension.
The symbols and the calculus depend strongly on the chosen basis $\bbb$.
 Nevertheless, the basis is arbitrary and we shall see, moreover, 
 that the analogue of the Laplace operator 
does not depend on $\bbb$, under precise conditions.\\
Some points, which are important in pseudodifferential analysis, have still
to be solved. The construction has been completed in \cite {AJN}, but the 
 covariance has been proved only in the most simple case.
 Beals characterization has been treated in
 \cite{ALN}, whereas
the composition results have been obtained in a high (but finite)  dimensional
 setting  \cite{AN1}. 

\medskip

In the present article, we introduce another class of symbols, $S(Q_A)$,
 defined thanks to a quadratic form (Definition \ref{Classe-DD-def}). Indeed, 
the first classes could be used in quantum electrodynamics but only under a 
truncature  assumption (see \cite{ALN2}), which the  classes $S(Q_A)$ enable
us to lift. The new  classes $S(Q_A)$ do not depend on a basis either.
 We go further with the study 
of the classes  $S_m(\bbb,\eps)$, begun in  \cite{AJN}, where only the
 properties absolutely necessary in view of the construction
itself had been developped. 
In the same way, we prove comparable properties for the new classes. The aim is 
to define a semigroup of operators similar to the heat operator and to 
state the properties which will be needed to treat the 
composition of operators.

We first generalize a result of  \cite{AJN} about the existence of stochastic
 extensions for  the  classes $S_m(\bbb,\eps)$ (Proposition
 \ref{prop84-modifiee}) and prove that their symbols are Frechet-differentiable
 for sufficiently large  $m$. We then prove the existence
of stochastic extensions for symbols in the new classes $S(Q_A)$ (Proposition 
\ref{Classe-DD-ext}).
Next we study the Taylor's expansions for the symbols in both classes, 
in order to give expansions for the regular terms and for the rest. 
One states there stochastic extension properties for unbounded  functions
 (Proposition \ref{4.11}), for which, in certain cases, the Weyl calculus has 
been defined by other means in  \cite{AJN}. This will prove
indispensable to treat the polynomial terms in Taylor's expansions.\\
The next step is the construction of a semigroup similar to the heat semigroup,
 denoted by $(H_t)_{t\geq 0}$. For bounded Borel functions defined on the
 Wiener space $B$ itself, it is a classical notion,
 to which \cite{G-4} is
 almost entirely devoted and which is still being studied (\cite{HA}). It is 
 given by
$$
\forall x \in B,\quad H_t f(x)= \int_B f(x+y)\ d\mu_{B,t}(y)
$$
with our notations. When $f$ is bounded and uniformly continuous, $H_tf$
converges uniformly to $f$ when $t$ converges to $0$. If, moreover, $f$
is Lipschitz continuous on $B$, $H_tf$ has further
 differentiability properties.\\
But our construction requires this notion for symbols $f$ defined on the 
initial  Hilbert space $H$, which is slightly complicated 
 since $H$ is $\mu_{B,t}$-negligible in $B$.
 We then set:
$$
\forall x\in H,\quad H_t f(x)= \int_B \tilde{f}(x+y)\ d\mu_{B,t}(y),
$$ 
where $\tilde{f}$ is a stochastic extension of $f$ in a certain sense. For 
a completion of $H$ given by an arbitrary measurable norm, this function
 $\tilde{f}$ is not necessarily continuous or uniformly continuous.
 But we prove, for both
classes, the existence of a precise completion $B_A$ of $H$ in which the
 stochastic extensions 
have useful topological properties (Propositions \ref{BA-classeDD} and
\ref{BA-ancienneclasse}). 
 This allows us to use the classical theory of
\cite{G-4} and \cite{K}.
 The extension $B_A$ may be different from the extension
 $B$ initially chosen and is used temporary.
 Of course, one checks that the integral defining  $H_tf(x)$ does not depend
on the chosen Wiener space completing $H$.

The main results of this article are Theorems \ref{glavni-Sm}, \ref{pcpal-DD},
 which establish, for the classes   $S_m(\bbb,\eps)$ 
as for the  classes $S(Q_A)$, the existence of a Laplacian commuting with the
 heat operator and which is its infinitesimal operator.
This result  will play an important part in the composition of symbols.
Let us stress the following expansion for $f\in S(Q_A)$, 
$$
H_tf= f + \sum_{k=1}^N \frac{t^k}{k!} \left(\frac{1}{2}\Delta \right)^kf 
+ t^{N+1} R_N(t), 
$$
where the rest satisfies estimates independent of $t$. To conclude, we give, 
for the class $S(Q_A)$, an invariance property  (Proposition \ref{comp-phi})
 which will be useful to prove a covariance result.

\medskip
Section  \ref{Rappels} recalls the indispensable notions about Wiener space
 and Wiener measure. Then it gives the vital definitions and results about
 the Weyl calculus in an infinite dimensional setting. Section
 \ref{3-ext-stoch} recalls and states more precisely the results about
stochastic extensions for the classes $S_m(\bbb,\eps)$
 of \cite{AJN}. It proves similar results in the case of products of scalar
 products and for the classes  $S(Q_A)$, which are defined at this point.
It brings up the alternative definition of the Weyl calculus, as a quadratic
form, which enables us to use unbounded symbols.
In Section  \ref{4-Taylor} we prove the Frechet-differentiability of the 
symbols in the  classes  $S_m(\bbb,\eps)$ and we extend
stochastically the Taylor's expansions
of symbols of both classes. Section  \ref{5-chaleur} defines the  heat
 operator $H_t$ for functions initially defined on the Hilbert space 
(and which it is impossible to integrate on the Wiener space without
an extension). We establish the semigroup property  for both classes,
 together with useful properties of $H_t$ (infinitesimal generator, 
commutation). The technical results about classical integration
are stated in the appendix
(section \ref{Annexe}).

\smallskip
The author wants to thank L.Amour and J.Nourrigat for many
 fruitful discussions.

%%%%%%%%%%%%%%%%%%%%%%%%%%%%%%%%%%%%%%%%%%%%%%%%%%%%%%%%%%%%%%%%%
%%%%%%%%%%%%%%%%%%%%%%%%%%%%%%%%%%%%%%%%%%%%%%%%%%%%%%%%%%%%%%%%%
%%%%%%%%%%%%%%%%%%%%%%%%%%%%%%%%%%%%%%%%%%%%%%%%%%%%%%%%%%%%%%%%%

\section{The Weyl calculus on a Wiener space}\label{Rappels}

The construction of the Wiener space may be found in 
\cite{G-1,G-2,G-3,K}. The Weyl calculus on a Wiener space has been developed
 in \cite {AJN}. We juste recall here the notions which are necessary to read
 the present article.

The abstract Wiener space $(i,H,B)$ is a triple where $H$ is a real, separable,
infinite dimensional  Hilbert space, $B$ is a Banach space containing $H$ and
$i$ is the canonical injection (which is not always mentioned).
Moreover,  $H$ is continuously embedded in $B$ as a dense subspace. Sometimes,
$B$ itself is called the Wiener space, as opposed to $H$, when no confusion  
is possible. One denotes by $<\,\ >$ (or sometimes $\cdot$ ) and $|\ |$ 
the scalar product and the norm on $H$ and by 
$|| \ ||$ the norm on $B$.\\
One identifies $H$ with its dual space, so that $B' \subset H\subset B$,  each
 space being a dense subspace of the following one.
One denotes by  $\fff(X)$ the set of all finite dimensional subspaces
of a vector space $X$. If $E\in \fff(H)$, one denotes by $\pi_E$ 
the orthogonal projection of $H$ onto $E$.\\
It is impossible to extend to $H$ itself the Gaussian measure which is
 naturally defined on its finite dimensional subspaces. Nevertheless, 
if the norm   $|| \ || $ of $B$ has a property called measurability
 (see  Definition 4.4 Chap 1 \cite{K} or \cite{G-1}), one can construct
a  Gaussian  measure on the Borel $\sigma$-algebra of $B$.  
Let us denote by 
$$
d\mu_{\R^n,h} (x)= (2\pi h)^{-n/2} e^{-\frac{1}{2h} \sum_{i=1}^n x_i^2}
 \ d\lambda(x_1,\dots, x_n)
$$
the Gaussian measure with variance $h>0$ on $\R^n$.
A cylinder of $B$ is  a set of the form
\begin{equation}
\ccc= \{ x\in B : (y_1(x),\dots,y_n(x))\in A \},
\end{equation}
where $n$ is a positive integer,
$y_1,\dots, y_n$ are elements of $B'$ and $A$ is a Borel set of  $\R^n$.
One defines the measure of this cylinder setting 
\begin{equation}\label{mes-cyl}
\mu_{B,h}(\ccc)=\int_A d\mu_{\R^n,h}(x)
\end{equation}
in case the family $(y_1,\dots, y_n)$ is orthonormal with respect to the scalar
 product of $H$, which can always be assumed. The parameter $h$ represents
 the variance of the Gaussian measure and can also be considered as a
 semiclassical parameter in the Weyl calculus.
One can prove that this measures extends as a probability measure, still
 denoted by  $\mu_{B,h}$, on the $\sigma$-algebra generated by the cylinders
 of $B$,  which is the Borel $\sigma$-algebra of $B$. The same definition,
but starting from cylinders of $H$, yields a pseudomeasure which is not
$\sigma$-additive.

If $E\in \fff(B')$ has dimension $n$, one can identify $E$ and $ \R^n$ by
choosing a basis, orthonormal with respect to the scalar product of $H$ and 
thus define a measure $\mu_{E,h}$ on $E$. For every
function $\ph\in L^1(E, \mu_{E,h})$,  the transfer  theorem
\begin{equation}\label{7AJN}
\int_B \ph\circ P_E(x) d\mu_{B,h}(x) = \int_{E} \ph(u) d\mu_{E,h}(u).
\end{equation}

If $y$ is an element of  $B'$, it can be considered as a random variable on
   $B$. If $y$ is not zero one sees, using (\ref{mes-cyl}),
 that, for every Borel set 
$A$ of $\R$, 
$$
\mu_{B,h}(y\in A)=\int_Ae^{-\frac{v^2}{2h|y|^2}} (2\pi h |y|^2)^{-1/2} \ dv,
$$
which means that $y$ has the normal distribution
 $\nnn(0,\sigma^2= h|y|^2)$ \cite{K}. 
Up to the factor  $\sqrt{h}$, there exists an isometry from $(B',|\ |)$ in
$L^2(B,\mu_{B,h})$. It can be extended as an isometry from $H$
in $L^2(B,\mu_{B,h})$ and one denotes
 by   $\ell_a$ the image of an element $a$ of $H$. If $a\in B'$, $\ell_a=a$ is
a linear application but if $a\in H$, $\ell_a$ is only defined $\mu_{B,t}$-
almost everywhere and is not necessaryly linear. However, 
 $\ell_a(-x)= -\ell_a(x)$ and
 $\ell_a(x+y)= \ell_a(x) + a\cdot y$ for  $y\in H$.
\\
If $E\in \fff(H)$ has an orthonormal basis $(e_1\dots, e_n)$, one sets, for 
$x\in B$, 
\begin{equation}\label{(50AJN)}
\tilde{\pi}_E(x)= \sum_{j=1}^n \ell_{u_j}(x) u_j,
\end{equation}
 in keeping with the projection. Then, for all $a\in H$, 
\begin{equation}\label{(51AJN)}
a\cdot \tilde{\pi}_E(x)= \ell_{\pi_E(a)}(x).
\end{equation}
The functions $\ell_a$ satisfy the following identities, recalled in
\cite{AJN}.
If  $a = u+i v$, with $u$ and $v$ in $H$, then
\begin{equation}
 \int _B e^{ \ell_a (x)}  d\mu _{B,h} (x) = e^{ h \frac{a^2}{ 2}}.
\label{(43AJN)} 
\end{equation}
One has set $a^2 = |u|^2 - |v|^2 + 2i u\cdot v $.
 For all $a$ in $H$ and for all 
 $p\geq 1$:
\begin{equation} \int _B |\ell_a (x)|^p   d\mu _{B,h} (x) = 
 \frac{(2h)^{p/2}}{ \sqrt { \pi} } |a|^p \
\Gamma \left (\frac{p+1}{ 2} \right ). \label{(44AJN)}
\end{equation}
Setting 
\begin{equation}\label{K(p)}
 K(p)= 2^{1/2} \pi^{-1/2p}
\left( \Gamma \left(\frac{p+1}{ 2} \right)\right)^{1/p},
\end{equation}
one can write that
$
\Vert \ell_a \Vert_{L^p(B,\mu _{B,h})} =K(p) h^{1/2} |a|.
$
Notice that $K(2)=1$.
One sees, too, that for all $a$ and $b$ in $H$,
 \begin{equation}
\int_B e^{\ell_b(u)} |\ell _a (u)|^p d\mu_{B,h} (u) =
e^{h\frac{|b|^2}{ 2}} \int_{\R} |\sqrt {h} |a|v + h a\cdot b|^p d\mu_{\R,1}(v).
\label{(45AJN)}
 \end{equation}

Let us recall the theorem of Wick :
\begin{theo}\label{Wick} {\bf Wick}
 Let $u_1, ... u_{2p}$ be vectors of $H$ ($p\geq 1$). Let  $h>0$. Then
one has
\begin{equation}
 \int _B \ell _{u_1} (x) ... \ell _{u_{2p}} (x) d\mu _{B,h} (x) = h^p
\sum _{(\varphi, \psi)\in S_p} \prod _{j=1}^p < u_{\varphi (j)}, u_{\psi (j)}>
\label{(48AJN)}
\end{equation}
where $S_p$ is the set of all couples  $(\varphi, \psi)$ of injections from
$\{ 1,...,p\}$ into $\{1,...., 2p\}$ such that:
\begin{enumerate}
\item
 For all  $j\leq p$,  $\varphi (j) < \psi (j)$.
\item
 The sequence $(\varphi (j) ) _{(1\leq j \leq p)}$ is an increasing sequence.
\end{enumerate}
\end{theo}

The measure  $\mu_{B,h}$  transforms, under translation of a vector $a$
 belonging to $H$, into another measure which is absolutely continuous with
 respect to the former one.
 More precisely,
 for all  $g\in L^1(B, \mu_{B,h}) $, one has, for all $a$
in $H$:
\begin{equation}
\int_B g(x) d\mu_{B,h} (x)= e^{-\frac{1}{2h}|a|^2}
\int_B g(x+a) e^{ - \frac{1}{h}\ell_a (x) }  d\mu _{B,h}(x).\label{(46AJN)} 
\end{equation}
But if the translation vector $a$ belongs to $B$, or if the variance parameter
 $h$ changes, both measures are mutually singular.\\
The Weyl calculus on the Wiener space has been constructed in two different 
ways. One of the constructions is rather similar to the classical definition,
 in that 
it relies on classes of symbols which satisfy differentiability conditions
and it yields operators which are bounded on a $L^2$ space. We will work in
 this frame  most of the time. We do not have, though, an integral definition
 of $Op(f)u$,
neither on $H$ nor on $B$.
The symbols are functions defined on $H^2$ by  Definition \ref{1.3AJN}. It is
 possible  - and necessary - to extend them to functions defined on $B^2$
 according to the definition below. This notion is classical in the theory
 of Wiener spaces
 (see \cite{G-1,G-2,G-3}, \cite{RA}, \cite{K}).
\begin{defi}\label{4.4AJN} 
 Let $(i, H,B)$  be an abstract Wiener space such that $B' \subset H\subset B$
 (The inclusion $i$ will be omitted). Let $h$ be a positive real number.
\begin{enumerate}
\item
 A  Borel function $ f$, defined on $H$, admits a stochastic
 extension $\widetilde f$ with respect to the measure  $\mu_{B,h}$ if, for
every increasing sequence $(E_n)$ in ${\cal F}(H)$, whose union is dense in $H$,
the sequence of functions $ f \circ \widetilde \pi _{E_n}$ (where 
 $\widetilde \pi _{E_n}$ is defined by  (\ref{(50AJN)})) 
 converges in
probability with respect to the measure $\mu_{B,h}$  to
$\widetilde  f$.  In other words, if, for every  $\delta>0$,
\begin{equation}
 \lim _{n\rightarrow + \infty } \mu_{B,h} \left(  \left \{
x \in B,\ \ \ \  | f \circ \widetilde \pi_{E_n } (x) -  \widetilde f (x) |
> \delta \right \} \right ) = 0. \label{(52AJN)}
\end{equation}
\item 
A function $f$ admits a stochastic extension 
 $\widetilde f$ in  $L^p (B,\mu _{B,h})$ 
($1\leq p< \infty$) if, for every increasing sequence $(E_n)$ in ${\cal F}(H)$,
 whose union is dense in $H$, the functions  $f \circ\widetilde \pi_{E_n}$
are in  $L^p (B,\mu_{B ,h})$ and if the sequence  $f \circ \widetilde \pi_{E_n}$
converges in $L^p (B,\mu _{B,h})$ to $\widetilde f$.
\end{enumerate}
 One defines likewise 
the stochastic extension of a function on $H^2$ to a function on  $B^2$.
\end{defi}
One can check, for example thanks to (\ref{(44AJN)}), that $\ell_a$ is the
 stochastic extension of the scalar product with $a$ and that $\tilde{\pi}_E$ 
is the stochastic extension of  $\pi_E$ in $L^p$.
The stochastic extension can be obtained in a more topological manner 
(see \cite{K}, chap. 1, par. 6). Let us draw attention to a result 
about extensions of holomorphic functions  (Theorem 8.8 \cite{AJN}),
obtained by martingale methods and proving a property announced by 
  \cite{K-R}.

The symbol classes used in   \cite{AJN} share derivability properties and
estimates with the classes 
of the Calder\'on-Vaillancourt Theorem:
\begin{defi}\label{1.3AJN}
 Let  $(i, H,B)$  be an abstract Wiener space
such that $B' \subset H\subset B$. Let $\bbb= (e_j )_{(j\in \Gamma)}$ 
be a  Hilbert basis of $H$, each vector belonging to
$B'$, indexed by a countable set $\Gamma$. Set $u_j = (e_j,0)$ and
$v_j = (0, e_j)$ $(j\in \Gamma)$. A multi-index is a map
 $(\alpha,\beta )$ from $\Gamma $ into
$\N \times \N$ such that $\alpha_j = \beta _j = 0$ except for a
finite number
 of indices. Let $M$ be a nonnegative real number, $m$  a nonnegative integer
 and  $\varepsilon = (\varepsilon_j )_{(j \in \Gamma)}$ a family of
nonnegative real numbers.  One denotes by $ S_m(\bbb, M, \varepsilon)$ the
set of bounded continuous functions $ F:H^2\rightarrow {\bf
C}$ satisfying the following condition. For every multi-index
 $(\alpha,\beta)$ of depth $m$, that is to say such that $0 \leq \alpha_j \leq
m$ and $0 \leq \beta_j \leq m$ for all $j\in \Gamma$, the following
derivative
\begin{equation}
\partial_u^{\alpha}\partial_v^{\beta}  F =  \left [\prod _{j\in \Gamma }
\partial _{u_j} ^{\alpha_j} \partial _{v_j} ^{\beta_j}\right ]  F  \label{(1.15)}
\end{equation}
is well defined, continuous on
 $H^2$ and satisfies, for every $(x,\xi)$ in  $H^2$
\begin{equation}
\left | \left [\prod _{j\in \Gamma } \partial_{u_j} ^{\alpha_j}
 \partial _{v_j} ^{\beta_j}\right ]  F(x,\xi)
\right |    \leq M \prod _{j\in \Gamma } \varepsilon_j ^{\alpha_j +
\beta_j}\ . \label{(1.16)} 
\end{equation}
\end{defi}

One recalls the following very useful property, stated in the
 proof of Proposition  4.14  of \cite{AJN}. If $\eps$ is square summmable,
   every  function $ F$ in $ S_1(\bbb,M, \varepsilon)$
verifies, for all $X$ and $V$ in $H^2$, a Lipschitz condition:
\begin{equation}\label{LIP}
| F(X+ V) -  F(X)| \leq M |V|\sqrt {2}  \left [
\sum_{j\in \Gamma} \varepsilon_j^2 \right ]^{1/2}.
\end{equation}

It is more convenient to represent classes of symbols as  vector spaces.
\begin{defi}\label{ancienne-classe}
Let  $\varepsilon $ be a sequence of positive real numbers and let $m\in \N$.
 One sets $S_m(\bbb,\varepsilon)= \bigcup_{M\geq 0} S_m(\bbb, M,\varepsilon)$.
For $F\in S_m(\bbb,\varepsilon)$ one sets
$||F||_{m,\varepsilon}= \inf\{ M\geq 0 : F \in  S_m(\bbb, M,\varepsilon)\}$.
\end{defi}
Remark that  $S_m(\bbb,\varepsilon)$, equipped with $||\ ||_{m,\varepsilon}$,
 is a Banach space. Setting 
 $S^{\infty}({\cal B},\varepsilon)= \bigcap_{m=0}^{\infty} S_m(\bbb,\varepsilon)$,
one can, classically, define a distance by
$
\displaystyle d(F,G)= \sum_{m=0}^{\infty} 2^{-m} 
\frac{ || F-G||_{m,\varepsilon}}{1+|| F-G||_{m,\varepsilon}}.
$
Then  $(S^{\infty}({\cal B},\varepsilon),d)$ is complete.\\

An alternative construction of the Weyl calculus uses an analogue of the Wigner
function in order to associate a quadratic form with a function  $\tilde{F}$
defined, this time, on $B^2$. This quadratic form is applied to cylindrical 
functions, depending on a finite number of variables. Let us only recall 
that this construction requires of $\tilde{F}$ to belong to 
 $L^1 ( B^2, \mu_{B^2 ,h/2})$ and to be such that there exists a 
nonnegative integer $m$ such that
\begin{equation}
 N_m ( \widetilde{F})=  \sup _{Y\in H^2 }
\frac{\Vert \widetilde{F}(\cdot + Y)\Vert_{L^1( B ^2,\mu_{B^2,h/2})}} {(1+|Y|)^m }
<+\infty.
  \label{(1.12AJN)} 
\end{equation}
This norm is finite if the function $\tilde{F}$ is bounded or if 
it is a polynomial
 expression of degree $m$ with respect to functions
 $(x,\xi) \rightarrow \ell _a (x) + \ell _b (\xi)$, with  $a$ and $b$
in $H$, as we shall see in  Subsection \ref{ext-prod-scal}.

These approaches complement one another. The most classical enables us to work 
on $L^2$ spaces on $B$, but the symbol has to be bounded, the other 
one allows us to use non bounded symbols, but the domain of the quadratic forms
contains only cylindrical functions. Both definitions coincide under 
certain conditions (Theorem 1.4 \cite{AJN}).

%%%%%%%%%%%%%%%%%%%%%%%%%%%%%%%%%%%%%%%%%%%%%%%%%%%%%%%%%%%%%%%%%%%%%%%%
%%%%%%%%%%%%%%%%%%%%%%%%%%%%%%%%%%%%%%%%%%%%%%%%%%%%%%%%%%%%%%%%%%%%%%%%
%%%%%%%%%%%%%%%%%%%%%%%%%%%%%%%%%%%%%%%%%%%%%%%%%%%%%%%%%%%%%%%%%%%%%%%%
%%%%%%%%%%%%%%%%%%%%%%%%%%%%%%%%%%%%%%%%%%%%%%%%%%%%%%%%%%%%%%%%%%%%%%%%

\section{Stochastic extensions }\label{3-ext-stoch} %3

%%%%%%%%%%%%%%%%%%%%%%%%%%%%%%%%%%%%%%%%%%%%%%%%%%%%%%%%%%%%%%%%%%%%%%%%
%%%%%%%%%%%%%%%%%%%%%%%%%%%%%%%%%%%%%%%%%%%%%%%%%%%%%%%%%%%%%%%%%%%%%%%%
%%%%%%%%%%%%%%%%%%%%%%%%%%%%%%%%%%%%%%%%%%%%%%%%%%%%%%%%%%%%%%%%%%%%%%%%

\subsection{Stochastic extensions of symbols in  $S_m(\bbb,\eps)$}%3.1

We first generalize a proposition  stated in \cite{AJN} (Proposition 8.4)
 in the case when $p=1$.
\begin{prop}\label{prop84-modifiee}
Let $F$ be a function in $S_1(\bbb,\varepsilon)$, with respect to a
 Hilbert basis
$\bbb= (e_j)_{(j\in \Gamma )}$, where the sequence $(\varepsilon_j) _{(j\in \Gamma)}$ 
is summable. Then, for every positive $h$ and every $q\in [1,+\infty[$,
 $F$ admits a stochastic extension in $L^q(B^2,\mu _{B^2,h})$.\\
 Moreover, for all $h_0>0$ and $q_0\in ]1,+\infty[$,  
there exists a function $\widetilde{F}$ which is the stochastic extension
of $F$ in  $L^q(B^2,\mu _{B^2,h})$ for all
$h\in ]0,h_0]$ and $q\in [1,q_0]$.\\
For any $E\in \fff(H^2)$, we then have the inequality :
$\forall (h,q)\in   ]0,h_0]\times  [1,q_0]$,
\begin{equation}
||F\circ \tilde{\pi}_E-\widetilde{F}||_{L^q(B^2,\mu_{B^2,h})} 
\leq  ||F||_{1,\eps} K(q)h^{1/2}
 \sum _{j=1}^{\infty} \varepsilon _j
 \Big (|u_j-\pi_{E}(u_j)| +|v_j-\pi_{E}(v_j) |\Big ). \label{ineg84}
\end{equation}
\end{prop}

{\it  Proof.} 
 Let $(E_n)$ be an increasing sequence of
 ${\cal F}(H^2)$, whose union is dense in $H^2$. For all $m$ and $n$ such
that $m<n$, let $S_{mn}$  be the orthogonal supplement of  $E_m$ in  $E_n$.
We can state an inequality analogous to  the inequality (120) of \cite{AJN}:
\begin{equation}\label{120.q}
|| F\circ \tilde{\pi}_{E_m}-F\circ \tilde{\pi}_{E_n} ||_{L^q(B^2,\mu_{B^2,h})}
\leq ||F||_{1,\eps} K(q)h^{1/2}
 \sum _{j=1}^{\infty} \varepsilon _j
 \Big (|\pi_{S_{mn}}(u_j)| +|\pi_{S_{mn}}(v_j) |\Big ).
\end{equation}

Indeed, one just needs to replace  $L^1$ by $L^q$ in the original proof, 
since the only changes take place in the explicit $L^q$ norms of the 
$\ell_a$ functions appearing there.
This inequality proves that $F\circ \widetilde \pi_{E_m}$ is a Cauchy sequence
 in $L^q(B^2,\mu_{B^2,h})$ and one can verify that the limit does not depend on
 the sequence $(E_n)$.  \\

To construct a representant of the stochastic extension common to all 
 $(h,q)\in ]0,h_0]\times [1,q_0]$, let us assume that $(E_n)$ is an increasing
 sequence of elements of $\fff(B')$.
The right term of (\ref{120.q}) is smaller than an expression $C(m,n)$ which
 depends, neither on $h$, nor on $q$. This allows us to construct 
an increasing sequence $(n_i)_i$ satisfying $C(n_{i+1},n_i)<2^{-i-1}$ and a 
sequence of functions $(F_N)_N$ defined by
$$
F_N:= F\circ \tilde{\pi}_{E_{n_{1}}}+
\sum_{j=1}^{N} \left( F\circ \tilde{\pi}_{E_{n_{j+1}}}- 
F\circ \tilde{\pi}_{E_{n_{j}}} \right)\ {\rm on }\ B^2,
$$
exactly as in the classical proof of the Riesz-Fisher Theorem.
The functions $F_N$ are defined everywhere on $B^2$ and independent 
of $(h,q)$. The limit  $\tilde{F}$ of this sequence is the representant we are 
looking for and it takes finite values on a subset of $B^2$ whose 
  $\mu_{B^2,h}$-measure is $1$  for all  $h\leq h_0$.\\
The last inequality  is a consequence of (\ref{120.q}), with  
 $E_m=E$ and letting $n$ grow to infinity.
\hfill $\square$

\begin{coro}\label{cor-prop84-modifiee}
Let $h_0$ be a positive real number. 
There exists a function $\widetilde{F}$  which is the stochastic extension
 of $F$ in $L^q(B^2,\mu_{B^2, h})$ for all $h\in ]0,h_0]$ and all
  $q\geq 1$. Inequality  (\ref{ineg84}) still holds.  
\end{coro}
{\it  Proof.}
Denote by  $\widetilde{F}_2$ (resp.  $\widetilde{F}_n$)
the stochastic expansion given by Proposition \ref{prop84-modifiee}, for
$h\in ]0,h_0]$ and  $q\in[1,2]$ (resp. $q\in[1,n]$). Let 
 $(E_s)$ be an increasing sequence of $\fff(H^2)$,  whose union is dense in 
 $H^2$.
One then has, for $q\leq 2$,
$$
\lim_{s\rightarrow \infty}
 || F\circ \tilde{\pi}_{E_s}- \widetilde{F}_2  ||_{q,h} =0,\quad
\lim_{s\rightarrow \infty}
 || F\circ \tilde{\pi}_{E_s}- \widetilde{F}_n  ||_{q,h} =0.
$$

Consequently $\widetilde{F}_2=\widetilde{F}_n$ $\mu_{B^2,h}$ almost everywhere.
 It follows that  $\widetilde{F}_2\in  L^q(B^2,\mu_{B^2, h})$ for 
$q\leq n$  and that the convergence is true. To obtain the inequality
one similarly replaces   $\widetilde{F}_n$
by $\widetilde{F}_2$.\hfill $\square$

Let us state another consequence of the proof of Proposition 8.4 of \cite{AJN}:

\begin{coro}\label{op-prol}
If $F\in S_1(\bbb,\varepsilon)$ where $\varepsilon $ is summable, 
then for all  $h>0$ and all $p\in [1,+\infty[$,
$$
| \widetilde{F}|\leq  ||F||_{1,\eps}\quad \mu_{B^2,h}-{\rm a.s.}\ , \qquad
 ||\widetilde{F}||_{L^p(B^2,\mu_{B^2,h})}\leq ||F||_{1,\eps}.
$$
Let us denote by  $\ppp$ the operator associating, with a function in
 $S_1(\bbb,\varepsilon)$, its stochastic extension in  $L^p(B^2, \mu_{B^2,t})$.
 This operator is thus linear, bounded and its norm is 
smaller than $1$.
\end{coro}
{\it  Proof.}
For every increasing sequence  $(E_n)_n$ of $\fff(H^2)$, whose union is dense 
in  $H^2$, the sequence  $(F\circ \tilde{\pi}_{E_n})_n$ converges to
$\widetilde{F}$ in $L^p$.  Since $|F\circ \tilde{\pi}_{E_n}|$  is smaller than 
 $ ||F||_{1,\eps}$ $\mu_{B^2,h}$ almost everywhere on $B^2$ (on the domain where
  $ \tilde{\pi}_{E_n}$ is defined or on  $B^2$ if
$E_n\subset B'$), so is $| \widetilde{F}|$.
Moreover, $||F\circ \tilde{\pi}_{E_n}||_{L^p(B^2,\mu_{B^2,h})}\leq  ||F||_{1,\eps}$ 
and letting $n$ grow to infinity yields  
 $|| \widetilde{F}||_{L^p(B^2,\mu_{B^2,h})}\leq  ||F||_{1,\eps}$.
\hfill $\square$

\begin{coro}\label{ext-transl}
If $F\in S_m(\bbb,\varepsilon)$ with $m\geq 1$ and $\varepsilon$ summable, if
$\widetilde{F}$ is the stochastic extension in the  $L^p$ given above, if
 $Y\in H^2$, then $\tau_YF$ admits $\tau_Y \widetilde{F}$  as a stochastic
 extension in the $L^p(B^2,\mu_{B^2,h})$ for $h>0$ and  $p\in [1,+\infty[$.
\end{coro}

{\it  Proof.}
Let $(E_j)$  be an increasing sequence  of $\fff(H^2)$, whose union is dense
 in  $H^2$. If we denote by an index $p$ the 
 $L^p(B^2,\mu_{B^2,h})$ norm, we obtain that
$$
|| \tau_Y \widetilde{F} -(\tau_Y F)\circ \widetilde{\pi}_{E_j} ||_p
\leq
|| \tau_Y \widetilde{F} -\tau_Y ( F\circ \widetilde{\pi}_{E_j}) ||_p+
||\tau_Y (F\circ\widetilde{\pi}_{E_j}) -(\tau_Y F)\circ\widetilde{\pi}_{E_j}||_p.
$$
For all  $p'>p$, the inequality
$$
|| \tau_Y \widetilde{F} -\tau_Y ( F\circ \widetilde{\pi}_{E_j}) ||_p 
= \left( \int_{B^2} |\widetilde{F} - F\circ \widetilde{\pi}_{E_j}|^p(X)
 e^{\frac{1}{h}\ell_Y(X)} d\mu_{B^2,h}(X) e^{-\frac{1}{2h} |Y|^2} \right)^{1/p}
\leq || \widetilde{F} - F\circ \widetilde{\pi}_{E_j}||_{p'}
 e^{\frac{|Y|^2}{2h(p'-p)}}
$$
holds true, thanks to the translation change of variables  (\ref{(46AJN)}),
 to H\"older's inequality and to the formula (\ref{(43AJN)}).
The second term converges to $0$ as well, thanks to (\ref{LIP}),
which enables us to give an upper bound of  
 $$
|F(\tilde{\pi}_{E_j}(X+Y)) - F(\tilde{\pi}_{E_j}(X)+Y)|=
|F(\tilde{\pi}_{E_j}(X)+\pi_{E_j}(Y)) - F(\tilde{\pi}_{E_j}(X)+Y)|.
 $$
\hfill $\square$ 

\begin{rem}\label{rem-ext-trans}
The result above holds for every globally Lipschitz continuous function $F$ 
admitting  a stochastic extension $\widetilde{F}$
in $L^p$ for every $p\in [1,+\infty[$ and $h>0$.
\end{rem}

%%%%%%%%%%%%%%%%%%%%%%%%%%%%%%%%%%%%%%%%%%%%%%%%%%%%
%%%%%%%%%%%%%%%%%%%%%%%%%%%%%%%%%%%%%%%%%%%%%%%%%%%%%%%%%%%%%%%%%%%%%%%%%

\subsection{Symbol classes defined thanks to a quadratic form} % 3.2

\begin{defi}\label{Classe-DD-def}
Let $A$ be a linear, selfadjoint, nonnegative, trace class application on a 
Hilbert space $H$. For all $x\in H$ one sets $Q_A(x)=  < Ax , x > $.
Let $S(Q_A) $ be the class of all functions $f\in C^{\infty } (H)$ such that 
there exists $C(f) >0$ satisfying:
\begin{equation}\label{Classe-DD-def-ineq}
\begin{array}{lll}
 \forall x\in H,\ 
|f  (x)  | \leq C (f),\\
\displaystyle
\forall m\in\N^*,\ \forall x\in H,\  \forall (U_1,\dots,U_m)\in H^m,\ 
|(d^m f ) (x) (U_1 , ... , U_m ) | \leq C (f) \prod _{j=1}^m
 Q_A( U_j) ^{\frac{1}{2}}.
\end{array}
\end{equation}
The smallest constant $C(f)$ such that  (\ref{Classe-DD-def-ineq})
holds is denoted by $\Vert f \Vert_{Q_A}$. 
\end{defi}

Notice that  $S(Q_A)$, equipped with the norm  $||\ ||_{Q_A}$, is a Banach space.
One can also check that, if $A$ and $B$ satisfy the conditions of
Definition \ref{Classe-DD-def}, their product  belongs to  $S(Q_{2(A+B)})$ with 
$$
|| fg||_{Q_{2(A+B)}} \leq || f||_{Q_{A}}|| g||_{Q_{B}}.
$$
Moreover, if $A$ is as above but defined on $H^2$, the class $S(Q_A)$ 
is included in a class $S_{\infty}(\bbb,\eps)$  for any orthonormal basis 
 $\bbb= (e_j)$ of  $H$, with $\eps_j= \max( Q_A(e_j,0)^{1/2}, Q_A(0,e_j)^{1/2})$.
Since the sequence $\eps$ is only square summable, the existence results for
 the stochastic extensions must be obtained otherwise.
\begin{rem} \label{avenant}
The constant $C(f)$ in the Definition \ref{Classe-DD-def} could depend on the 
 order $m$. Some results are still valid with a less restrictive class
 satisfying
$$\forall m\in\N^*,\ \exists C_m= C_m(f) :\ 
 \forall x\in H,\  \forall (U_1,\dots,U_m)\in H^m,\ 
|(d^m f ) (x) (U_1 , ... , U_m ) | \leq C_m (f) \prod _{j=1}^m
 Q_A( U_j) ^{\frac{1}{2}}.
$$
\end{rem}

\begin{lemm}\label{ineg-Q}
For $E\in \fff(H)$ and  $h>0$,
$y\mapsto Q_A(\tilde{\pi}_E(y))^{\frac{1}{2}}$ belongs to $L^p(B,\mu_{B,h})$ 
for all $p\in[1,+\infty[$. More precisely, if  $(u_j)$ is a Hilbert basis of
 $H$ whose vectors are eigenvectors of $A$ (or belong to ${\rm Ker}(A)$)
 and if one denotes by
 $\lambda _j $ the corresponding eigenvalues, one obtains
$$
||Q_A^{\frac{1}{2}}\circ\tilde{\pi}_E||_{L^p(B,\mu_{B,h})} \leq C(p)
\left( \sum_{0}^{\infty} \lambda_j |\pi_E(u_j)|^{\alpha(p)} \right)^{1/{\alpha(p)}}
 h^{\frac{1}{2}},
$$
with
\begin{equation}\label{constantesQ}
\begin{array}{llll}
\displaystyle C(p)= K(p)\left( \sum_{0}^{\infty} \lambda_j \right)^{\frac{1}{2}
 -\frac{1}{p}}  & \alpha(p)= p  & {\rm for} \ p> 2\\
\displaystyle C(p)= 1  & \alpha(p)= 2  & {\rm for} \ p\leq 2,\\
\end{array}
\end{equation}
the constant  $K(p)$ being defined by (\ref{K(p)}).
\end{lemm}

{\it Proof.}
By decomposing $A$ on its eigenvector basis, one obtains that
$$
Q_A(\tilde{\pi}_E(y))=
\sum_{j=0}^{\infty} \lambda_j (u_j\cdot  \tilde{\pi}_E(y))^2
=\sum_{j=0}^{\infty} \lambda_j (\ell_{\pi_E(u_j)})^2,
$$
using (\ref{(51AJN)}). For $p=2$ it suffices to integrate this equality and to
 use  (\ref{(44AJN)}).
For $p>2$, one uses Jensen's inequality  for a probability measure on $\N$.
Set $S=  \sum_{0}^{\infty} \lambda_j$. One then has 
$$
Q_A(\tilde{\pi}_E(y))^{\frac{p}{2}}
=\left( \sum_{j=0}^{\infty} \frac{\lambda_j}{S} \ 
S (\ell_{\pi_E(u_j)})^2\right)^{\frac{p}{2}}
\\
\leq
 \sum_{j=0}^{\infty} \frac{\lambda_j}{S} \ 
( S (\ell_{\pi_E(u_j)})^2) ^{\frac{p}{2}},
$$
and it remains to integrate. Finally, for $p\in [1,2[$, one applies  H\"older's
 inequality.
\hfill $\square$\\

\begin{rem}\label{rem-maj-Q}
One can give an upper bound for
$|| Q_A\circ \tilde{\pi}_E||_{L^p(B,\mu_{B,h})}$, which does not depend on $E$:
$$
||Q_A^{\frac{1}{2}}\circ\tilde{\pi}_E||_{L^p(B,\mu_{B,h})} \leq C(p)
\left( \sum_{0}^{\infty} \lambda_j \right)^{1/{\alpha(p)}}
 h^{\frac{1}{2}}.
$$
\end{rem}

One can prove the following result.
\begin{prop}\label{Classe-DD-ext}
Let $h>0$ and let $p\in[1,+\infty[$.
Every function $f$ belonging to  $S(Q_A)$ admits a stochastic extension 
 $\widetilde f $ in $L^p(B,\mu_{B,h})$. The function $\widetilde f $ is bounded
 $\mu_{B,h}$ almost everywhere by $||f||_{Q_A}$. \\
Moreover, for all $E\in \fff(H)$,
\begin{equation}\label{Classe-DD-ext-ineg}
  ||f \circ \tilde{\pi}_{E} - \widetilde{f} (x) ||_{L^p(B, \mu _{B , h })} 
\leq
C(p) h^{\frac{1}{2}}  \Vert f \Vert _{Q_A} 
\left(\sum_{j\geq 0}\lambda _j |\pi _{E} (u_j )-u_j|^{\alpha(p)} \right)^{1/\alpha(p)}, 
\end{equation}
with the notations of Lemma \ref{ineg-Q}.
\end{prop}

 {\it Proof.} Let  $(E_n)$ be an increasing sequence  of $\fff(H)$, whose 
union is dense in $H$.
 Let  $f$ be in $S(Q_A)$. Let  $m$ and $n$ be such that $m<n$. Let $S_{mn}$ be
 an orthogonal supplement of  $E_m$ in  $E_n $. Then
 $$ f (\tilde{\pi}_{E_n } (x) ) - f (\tilde{\pi}_{E_m} (x) ) = \int _0^1
 (df) ( \tilde{\pi}_{E_m} (x) + \theta \tilde{\pi}_{S_{mn} } (x) )
 ( \tilde{\pi}_{S_{mn} } (x) ) d\theta .$$  
 Hence 
 $$
 |f (\tilde{\pi}_{E_n } (x) ) - f (\tilde{\pi}_{E_m} (x) )|
 \leq
|| f||_{Q_A} \int_0^1  Q_A ( \tilde{\pi}_{S_{mn} } (x) )^{\frac{1}{2}}  d\theta
=|| f||_{Q_A}   Q_A ( \tilde{\pi}_{S_{mn} } (x) )^{\frac{1}{2}}.
 $$   
This implies that 
$$
||f\circ \tilde{\pi}_{E_n }-f\circ \tilde{\pi}_{E_m }||_{L^p(B,\mu_{B,h})}
\leq  \Vert f \Vert _{Q_A} 
||   Q_A^{\frac{1}{2}} \circ \tilde{\pi}_{S_{mn} }  ||_{L^p(B,\mu_{B,h})}.$$
Using the preceding Lemma \ref{ineg-Q}, one gets that
$$
||f\circ \tilde{\pi}_{E_n }-f\circ \tilde{\pi}_{E_m }||_{L^p(B,\mu_{B,h})}
\leq C(p) h^{\frac{1}{2}}  \Vert f \Vert _{Q_A} 
\left(\sum_{j\geq 0}\lambda _j |\pi _{S _{mn}} (u_j )|^{\alpha(p)} \right)^{1/\alpha(p)}.
$$
The right term converges to $0$ when $m$ grows to infinity, according to
 the dominated convergence Theorem. Indeed, for all $j$,
$ |\pi _{S _{mn}} (u_j ) |$ converges to $0$ when $m$ grows to infinity,
$ |\pi _{S _{mn}} (u_j ) |^{\alpha(p)}  \leq 1$ and the series $\sum\lambda _j$
 converges. The sequence $(f (\tilde{\pi}_{E_n }))_n$ is therefore a Cauchy
 sequence in $L^p(B,\mu _{B , h } )$. One can verify that its limit, in
 $L^p(B,\mu _{B , h } )$, does not depend on the sequence $(E_n)$. 
Since the function $|f \circ\tilde{\pi}_{E_n }|$ is almost everywhere smaller
 than $||f||_{Q_A}$, so is its limit. Finally, taking $E=E_m$  in one of
 the above inequalities and letting $n$ converge to infinity yields 
(\ref{Classe-DD-ext-ineg}).\hfill $\square$

\begin{rem}
This result holds true for the class of Remark \ref{avenant}, with
 $\max( C_0(f), C_1(f))$ instead of  $||f||_{Q_A}$  in the estimates.
\end{rem}

\begin{prop}\label{prol-diff-DD}
Let $h>0$, $p\in[1,+\infty[$. Let $k$ be a positive integer and let  $x$ be a
 fixed point in $H$. Set $S= \sum\lambda_j$.
The function $y \mapsto d^kf(x)\cdot y^k$ defined on   $H$
 admits a stochastic extension in $L^p(B,\mu_{B,h})$.
Moreover, for all $E\in \fff(H)$, 
$$
 ||d^kf(x)\cdot \tilde{\pi}_E(y)^k- \ppp(y\mapsto d^kf(x)\cdot y^k)||_{p}
\leq 
k||f||_{Q_A}   C(pk)^k S^{ \frac{k-1}{\alpha(pk)} }h^{\frac{k}{2}} 
\left( \sum \lambda_s |{\pi}_E(u_s)- u_s |^{\alpha(pk)} 
\right)^{\frac{1}{\alpha(pk)}}.
$$
\end{prop}
{\it Proof.} 
Let $E, F\in\fff(H)$ with $E\subset F$. 
For all  $y\in B$, one has
$$
\begin{array}{lll}
\displaystyle
 d^kf(x)\cdot \tilde{\pi}_E(y)^k- d^kf(x)\cdot \tilde{\pi}_F(y)^k&
\displaystyle
= \sum_{j=1}^k d^kf(x)( \tilde{\pi}_E(y)^{j}, \tilde{\pi}_F(y)^{k-j}) -
d^kf(x)( \tilde{\pi}_E(y)^{j-1}, \tilde{\pi}_F(y)^{k-j+1}) \\
&\displaystyle
= \sum_{j=1}^k d^kf(x)( \tilde{\pi}_E(y)^{j-1},\tilde{\pi}_E(y)-
 \tilde{\pi}_F(y), \tilde{\pi}_F(y)^{k-j}).
\end{array}
$$
Using Definition \ref{Classe-DD-def}, one deduces that 
$$
\begin{array}{lll}
\displaystyle
 |d^kf(x)\cdot \tilde{\pi}_E(y)^k- d^kf(x)\cdot \tilde{\pi}_F(y)^k|
\leq &
\displaystyle
 \sum_{j=1}^k ||f||_{Q_A} Q_A^{\frac{j-1}{2}}(\tilde{\pi}_E(y) ) Q_A^{\frac{k-j}{2}}(\tilde{\pi}_F(y) )Q_A^{\frac{1}{2}}(\tilde{\pi}_E(y)-\tilde{\pi}_F(y)).
\end{array}
$$
Using H\"older's inequality  and  Remark \ref{rem-maj-Q}, one obtains that 
$$
\begin{array}{lll}
\displaystyle
 ||d^kf(x)\cdot \tilde{\pi}_E(y)^k- d^kf(x)\cdot \tilde{\pi}_F(y)^k||_{p}\\
\leq 
\displaystyle
\sum_{j=1}^k ||f||_{Q_A}  ||Q_A^{\frac{1}{2}}\circ \tilde{\pi}_E||_{pk}^{j-1}  
||Q_A^{\frac{1}{2}}\circ \tilde{\pi}_F||_{pk}^{k-j }  
||Q_A^{\frac{1}{2}}\circ (\tilde{\pi}_E-\tilde{\pi}_F)||_{pk}\\
\leq
\displaystyle
k||f||_{Q_A}   (C(pk)S^{\frac{1}{\alpha(pk)}}h^{\frac{1}{2}} )^{k-1}  C(pk) h^{\frac{1}{2}} 
\left( \sum \lambda_s |({\pi}_E-{\pi}_F)(u_s) |^{\alpha(pk)} 
\right)^{\frac{1}{\alpha(pk)}}. 
\end{array}
$$
Then one proceeds as in the preceding proposition, replacing $E$ and $F$ 
by the terms of an increasing sequence  of $\fff(H)$ whose union is dense in
 $H$ and whose first term is $E$.
\hfill $\square$

\begin{rem}  This result holds with the class defined by Remark \ref{avenant},
 with  $C_k(f)$ instead of $||f||_{Q_A}$.
\end{rem}

 A consequence of Lemma \ref{ineg-Q} is the following result, 
which partly generalizes  Proposition 8.7 of \cite{AJN}:
\begin{coro}
Let $h>0$ and $p\in [1,+\infty[$. The function $Q_A^{\frac{1}{2}}$ admits
 a stochastic extension in  $L^p(B,\mu_{B,h})$.
\end{coro}
{\it Proof.} As in the proof of Proposition \ref{Classe-DD-ext}, one
 introduces an increasing sequence  $(E_n)$ of $\fff(H)$. 
One denotes by $S_{mn}$ an orthogonal supplement of $E_m$ in $E_n$ if $m\leq n$.
Lemma \ref{ineg-Q} then  implies the inequality 
$$
||Q_A^{\frac{1}{2}}\circ\tilde{\pi}_{S_{mn}}||_{L^p(B,\mu_{B,h})} \leq C(p)
\left( \sum_{0}^{\infty} \lambda_j |\pi_{S_{mn}}(u_j)|^{\alpha(p)} \right)^{1/{\alpha(p)}}
 h^{\frac{1}{2}},
$$
which proves that  $(Q_A^{\frac{1}{2}}\circ\tilde{\pi}_{E_n})_n$ is a  Cauchy
sequence in $L^p(B,\mu_{B,h})$. \hfill $\square$

We finally may state the following result, which enables us to use another
Wiener space associated with $H$ than the space $B$ initially chosen.
\begin{prop}\label{BA-classeDD}
Let  $A$ be a linear, selfadjoint, nonnegative, trace class application in a
 Hilbert space $H$. There exists a measurable  norm  (see \cite{K} (Def.4.4) or
 \cite{G-1}), $|| \ ||_{A,n}$ on  $H$, and hence a completion $B_A$ of $H$ with respect to this norm, such that the following property is satisfied: if
$f$ belongs to the class $S(Q_A)$, then $f$ is uniformly continuous 
on $H$ with respect to the norm $|| \ ||_{A,n}$.\\
 The function $f$ admits a uniformly continuous  extension $f_A$ on $B_A$ and 
the stochastic extension $\tilde{f}$ of $f$  given by
 Proposition \ref{Classe-DD-ext}
is  equal to $f_A$ $\mu_{B,h}$- a.e.

\end{prop}
{\it Proof.}
If  $A$ is an injection, one sets  $||x||_{A,n}= <Ax,x>^{1/2}= Q_A(x)^{1/2}.$
If not, if $(e_n)_n$ is an orthonormal basis of ${\rm Ker}(A)$, one can add to 
$A$ the operator $C$ defined, for example, by 
$Cx= \sum_{n} e^{-n} <x,e_n> e_n$. The operator $A+C$ is  selfadjoint,
 nonnegative, trace class  and it is an injection. One then sets
$||x||_{A,n}= <(A+C)x,x>^{1/2}= Q_{A+C}(x)^{1/2}.$
It follows from Theorem 3 in  \cite{G-1} that $|| \ ||_{A,n}$ is a measurable
 seminorm. It is a norm since $A$ (or $A+C$) is injective.
Taylor's formula gives, in both cases, the inequality 
$$\label{lipschitz-DD}
|f(y)-f(x)| \leq ||f||_{Q_A} Q_A(y-x)^{1/2}  \leq  ||f||_{Q_A} ||x-y||_{A,n} ,
$$
which in turn implies the uniform continuity. 
The topological extension $f_A$ of $f$ is then uniformly continuous on $B_A$.
According to Theorem 6.3 (Chap 1 \cite {K}), $f_A$ and $\tilde{f}$
coincide almost everywhere.
\hfill $\square$

%%%%%%%%%%%%%%%%%%%%%%%%%%%%%%%%%%%%%%%%%%%%%%%%%%%%%%%%%%%%%%%%%%%%%%%%
%%%%%%%%%%%%%%%%%%%%%%%%%%%%%%%%%%%%%%%%%%%%%%%%%%%%%%%%%%%%%%%%%%%%%%%%
%%%%%%%%%%%%%%%%%%%%%%%%%%%%%%%%%%%%%%%%%%%%%%%%%%%%%%%%%%%%%%%%%%%%%%%%
%%%%%%%%%%%%%%%%%%%%%%%%%%%%%%%%%%%%%%%%%%%%%%%%%%%%%%%%%%%%%%%%%%%%%%%%
%%%%%%%%%%%%%%%%%%%%%%%%%%%%%%%%%%%%%%%%%%%%%%%%%%%%%%%%%%%%%%%%%%%%%%%%

\subsection{Scalar products and products of scalar products }% 3.3
\label{ext-prod-scal}

%%%%%%%%%%%%%%%%%%%%%%%%%%%%%%%%%%%%%%%%%%%%%%%%%%%%%%%%%%%%%%%%%%%%%%%%

\begin{lemm}
For all $a\in H$, the function $\varphi_a$ : $H\rightarrow \R$,
 $x\mapsto <x,a>$
admits the function $\ell_a$ as a stochastic extension in 
$L^p(B,\mu_{B,h})$, for all $p\in [1,\infty[$ and all $h>0$.
Moreover,  for all  $E\in \fff(H)$,
\begin{equation}\label{ext_prodscal}
|| \varphi_a \circ \widetilde{\pi_{E}} - \ell_a ||_{ L^p(B,\mu_{B,h})}
= K(p)h^{\frac{1}{2}} | \pi_{E}(a)-a|.
\end{equation}
\end{lemm}
{\it  Proof.}
Let $(E_j)_j$ be an increasing sequence  of
 $\fff(H)$ such that $\overline{\bigcup E_j} =H$. Since
$\varphi_a\circ \widetilde{\pi_{E_j}} = \ell_{\pi_{E_j}(a)}$, according to 
 (\ref{(51AJN)}), one obtains, using  (\ref{(44AJN)}), that, for a finite $p$ 
$$
|| \varphi_a \circ \widetilde{\pi_{E_j}} - \ell_a ||_{ L^p(B,\mu_{B,h})} =
||  \ell_{\pi_{E_j}(a)-a} ||_{ L^p(B,\mu_{B,h})}
=
K(p)h^{\frac{1}{2}} | \pi_{E_j}(a)-a|.
$$
This proves the convergence and the result.\hfill $\square$

We now study the products of scalar products. Let  $a_1,\dots, a_n$ be vectors
 of $H$. Let  $\alpha=(\alpha_1,\dots, \alpha_n)$   be a multiindex such that
 $\alpha_i > 0$ for every  $i$.
One defines the function $a^{\alpha} $ on  $H$ by
$$
a^{\alpha}(x)= \prod_{i=1}^n <a_i,x>^{\alpha_i}.
$$
\begin{prop}\label{polyn-prodscal}
For $h>0$ and $p\in [1,+\infty[$, the function $a^{\alpha}$ admits the function
$\prod_{i=1}^n \ell_{a_i}^{\alpha_i}$ as a stochastic extension in 
$L^p(B,\mu_{B,h})$.
Moreover,  for all  $E\in \fff(H)$, 
$$
\left\Vert  a^{\alpha}\circ \widetilde{\pi_{E}} -
\prod_{i=1}^n \ell_{a_i}^{\alpha_i}\right\Vert_{L^p(B,\mu_{B,h})}
\leq   K(p|\alpha|)^{|\alpha|} h^{|\alpha|/2} ( \max_{1\leq i\leq n}|a_i|)^{|\alpha|-1} \sum_{i=1}^n \alpha_i |\pi_E(a_i)-a_i|.
$$
\end{prop}

{\it Proof.}
One obtains the inequality by applying Lemma \ref{Hoelder-telescopique} stated
in the appendix to the $|\alpha|$ functions appearing in the products
 $a^{\alpha}\circ \tilde{\pi}_E$ and $\prod _{i=1}^n \ell_{a_i}^{\alpha_i}$
and by remarking that 
$$
|| \ell_{\pi_{E}(a_i)}||_{p|\alpha|}= K(p|\alpha|)h^{\frac{1}{2}} |\pi_{E}(a_i)|\leq 
K(p|\alpha|)h^{\frac{1}{2}} |a_i|.
$$
It then remains to replace $E$ by  an increasing sequence  of $\fff(H^2)$ 
 such that $\overline{\bigcup E_j} =H^2$ to obtain the
stochastic extension.\hfill $\square$\\

Besides, one can define a quadratic form associated with such a product
thanks to the following result:
\begin{prop}\label{NormeNm-fonctionsell}
Let $a_1,\dots, a_n,b_1,\dots,b_p$ belong to  $H$, let
$\alpha_1,\dots,\alpha_n,\beta_1,\dots ,\beta_p$ be positive integers and set
 $m=\max(\sum_1^n \alpha_j,\sum_1^p \beta_j)=\max(|\alpha|,|\beta|)$.
The function
$\displaystyle
\widetilde{F} : (x,\xi)\mapsto  \prod_{i=1}^n \ell_{a_i}^{\alpha_i}(x)
 \prod_{i=1}^p \ell_{b_i}^{\beta_i}(\xi) 
$
has a finite norm  $N_m$ defined by (\ref{(1.12AJN)}). More precisely, 
\begin{equation}
\begin{array}{cc}
\displaystyle
N_m(\widetilde{F})  = \sup_{Y\in H^2}
\frac{||\widetilde{F}(\cdot+Y)||_{L^1(B^2,\mu_{B^2,\frac{h}{2}})}}{(1+|Y|)^m} \\
\displaystyle  \leq \max(1,\sqrt{\frac{h}{2}})^{ |\alpha|+|\beta| }
\prod_1^n|a_j|^{\alpha_j} \prod_1^p |b_i|^{\beta_i}
\times \prod_1^n\left(  \int_{\R}(1+|v|)^{n\alpha_j} \
 d\mu_{\R,1}(v)\right)^{\frac{1}{n}} 
\prod_1^p \left(  \int_{\R}(1+|v|)^{p\beta_i} \ d\mu_{\R,1}(v)\right)^{\frac{1}{p}} 
\end{array}
\end{equation}
\end{prop}

{\it Proof.}
By the change of variables formula (\ref{(46AJN)}) one obtains
$$
\int_{B^2} |\widetilde{F}(X+Y)| d\mu_{B^2,\frac{h}{2}}(X)
\leq
e^{-\frac{1}{h}|y|^2}\int_B \prod_{j=1}^n|\ell_{a_j}(x)|^{\alpha_j}
e^{\frac{2}{h} \ell_y(x)} d\mu_{B,\frac{h}{2}}(x)
e^{-\frac{1}{h}|\eta|^2}\int_B \prod_{j=1}^p|\ell_{b_j}(\xi)|^{\beta_j}
e^{\frac{2}{h} \ell_{\eta}(\xi)} d\mu_{B,\frac{h}{2}}(\xi).
$$
H\" older's inequality yields
$$
A
:=e^{-\frac{1}{h}|y|^2}\int_B \prod_{j=1}^n|\ell_{a_j}(x)|^{\alpha_j}
e^{\frac{2}{h} \ell_y(x)} d\mu_{B,\frac{h}{2}}(x)
\leq 
e^{-\frac{1}{h}|y|^2} \prod_{j=1}^n\left( \int_B |\ell_{a_j}(x)|^{n\alpha_j}
e^{\frac{2}{h} \ell_y(x)} d\mu_{B,\frac{h}{2}}(x)
\right)^{1/n}.
$$
According to (\ref{(45AJN)}),
$$
A \leq e^{-\frac{1}{h}|y|^2} \prod_{j=1}^n\left(e^{ |y|^2/h}    \int_{\R}
| \sqrt{\frac{h}{2}} |a_j|v + < y, a_j >|^{n\alpha_j} \ d\mu_{\R,1}(v) 
\right)^{1/n}.
$$

One can factor  $|a_j| $ and, remarking that
$
\sqrt{\frac{h}{2}} |v| +|y| $ is smaller than
 $ \max(1, \sqrt{\frac{h}{2}}) (1+|v|)(1+|y|),
$
one gets 
$$
A\leq  \max(1,\sqrt{\frac{h}{2}})^{|\alpha| }\prod_1^n|a_j|^{\alpha_j}
(1+|y|)^{ |\alpha|}
\prod_1^n\left(  \int_{\R}(1+|v|)^{n\alpha_j} \ d\mu_{\R,1}(v)\right)^{1/n}.
$$
One treats the other factor similarly, which gives the desired result.
\hfill $\square$

%%%%%%%%%%%%%%%%%%%%%%%%%%%%%%%%%%%%%%%%%%%%%%%%%%%%%%%%%%%%%%%%%%%%%%%%%
%%%%%%%%%%%%%%%%%%%%%%%%%%%%%%%%%%%%%%%%%%%%%%%%%%%%%%%%%%%%%%%%%%%%%%%%%
%%%%%%%%%%%%%%%%%%%%%%%%%%%%%%%%%%%%%%%%%%%%%%%%%%%%%%%%%%%%%%%%%%%%%%%%%
%%%%%%%%%%%%%%%%%%%%%%%%%%%%%%%%%%%%%%%%%%%%%%%%%%%%%%%%%%%%%%%%%%%%%%%%%
\subsection{Stochastic extension in an integral} % 3.4
\begin{lemm}
Let  $f : [0,1] \rightarrow \R$  and  $g :[0,1]\times B^2\rightarrow \R $ 
be a measurable functions. Let $q\in [1,+\infty[$ and let $h>0$.
If
$$
\int_0^1 |f(s)| \left(\int_{B^2}|g(s,Y)|^q \ d\mu_{B^2,h}(Y) \right)^{1/q} \ ds
 <\infty,
$$
then  $Y\mapsto \int_0^1 f(s)g(s,Y)\ ds$ belongs to  $L^q(B^2,\mu_{B^2,h})$ and
$$
\left\Vert  \int_0^1 f(s)g(s,\cdot )\ ds\right\Vert_{L^q(B^2,\mu_{B^2,h})}\leq  
\int_0^1 |f(s)| \left(\int_{B^2}|g(s,Y)|^q \ d\mu_{B^2,h}(Y) \right)^{1/q} \ ds.
$$
\end{lemm}
{\it  Proof.}
When  $q=1$, the result is straightforward. If $q>1$, one introduces
a function of  $ L^{q'}(B^2,\mu_{B^2,h})$ (where $q'$ is the 
conjugate exponent of $q$) and one proves that the integral belongs to
 the dual space of $L^{q'}$. 
\hfill $\square$\\

\begin{prop}\label{prol-int}
Let $f : [0,1] \rightarrow \R$ be a continuous function, let 
 $G$ be in $ S_1(\bbb,\varepsilon)$, 
with  $\varepsilon$ summable and let us fix $X\in H^2$.
For all  $q\in [1,+\infty[$ and all $h\in]0,1]$, the function 
$$ \displaystyle
Y\mapsto \int_0^1 f(s) G(X+sY)\ ds,
$$
defined on  $H^2$, admits, as a stochastic extension in
 $L^q(B^2,\mu_{B^2,h})$, the function 
$$ \displaystyle
Y\mapsto \int_0^1 f(s) \widetilde{G}(X+sY)\ ds,
$$
defined on  $B^2$, where $\widetilde{G} $ is the stochastic extension of $G$ 
for all $L^r(B^2,\mu_{B^2,h}), (r,h)\in [1,+\infty[\times ]0,1]$.\\
Moreover, if $E\in \fff(H^2)$, one has the inequality
$$
\begin{array}{llll}
\displaystyle
\left \Vert
\int_0^1f(s)\left(G(X+s\widetilde{\pi_E}(\cdot)) - \widetilde{G}(X+s\cdot)\right)\ ds
  \right\Vert_{ L^q(B^2,\mu_{B^2,h})}\\
\displaystyle \leq
||G||_{1,\varepsilon} \left( \int_0^1|f(s)|\ ds\right) \left(\sqrt{2\sum_{\Gamma}\varepsilon_j^2 }
 |X-\pi_E(X)| + \sqrt{2h(q+2)} C e^{\frac{|X|^2}{2h}}
 \sum_{0}^{\infty} \varepsilon_j(|u_j-\pi_E(u_j)|+|v_j-\pi_E(v_j)| )
 \right),
\end{array}
$$
where the constant $C$ does not depend on the parameters. 
\end{prop}
% vérifié au moins deux trois fois, dont le 16 mai
{\it  Proof.}
One checks that all the functions
$(s,Y) \mapsto X+sY,  X+s\tilde{\pi}_E(Y)$ are measurable. Set 
$U_E=\int_0^1f(s)(G(X+s\widetilde{\pi_E}(\cdot)) -\widetilde{G}(X+s\cdot))\ ds$.
Using the lemma above, one sees that 
$$
\begin{array}{lll}
||U_E||_{L^q(B^2,\mu_{B^2,h})} &
\displaystyle
\leq \int_0^1 |f(s)|\  || G(X+s\widetilde{\pi_E}(\cdot)) - 
 G(\widetilde{\pi_E}(X+s \cdot)) ||_{L^q(B^2,\mu_{B^2,h})} \ ds \\
& + \displaystyle
 \int_0^1 |f(s)|\  ||  G(\widetilde{\pi_E}(X+s \cdot))-
\widetilde{G}(X+s \cdot) ||_{L^q(B^2,\mu_{B^2,h})} \ ds.\\
\end{array}
$$
Formula  (\ref{LIP}) proves that the first term is smaller than
$$
\int_{B^2}
|G(X+s\widetilde{\pi_E}(Y))-G(\widetilde{\pi_E}(X+s Y))|^q \ d\mu_{B^2,h} \leq 
\int_{B^2} \left(||G||_{1,\varepsilon} \sqrt{2\sum\varepsilon_j^2}
 |X-\pi_E(X)|\right)^q \ d\mu_{B^2,h}.
$$
For the second term, successive change of variables give
$$
\begin{array}{rrr}
\displaystyle
\int_{B^2}
|G(\widetilde{\pi_E}(X+s Y))-\widetilde{G}(X+sY) |^q \ d\mu_{B^2,h}(Y)
= 
\int_{B^2}
|G(\widetilde{\pi_E}(X+Z))-\widetilde{G}(X+Z) |^q \ d\mu_{B^2,s^2h}(Z)\\
= \displaystyle
\int_{B^2}
|G(\widetilde{\pi_E}(Z))-\widetilde{G}(Z) |^q e^{-\frac{|X|^2}{ 2s^2h}}
 e^{\frac{1}{s^2h} \ell_X(Z)}\
 d\mu_{B^2,s^2h}(Z).\\
\end{array}
$$
One then applies H\"older's inequality to the last term, raising 
$|G(\widetilde{\pi_E}(Z))-\widetilde{G}(Z) |^q$ to the power   $q'/q$ with  
$q'=q+ \frac{1}{s^2}$.This gives
$$
|| G(\widetilde{\pi_E}(X+s \cdot))-\widetilde{G}(X+s\cdot) ||_{L^q(B^2,\mu_{B^2,h})}
\leq e^{\frac{|X|^2}{2h}}
|| G\circ \widetilde{\pi_E}-\widetilde{G}||_{L^{q+\frac{1}{s^2}}(B^2,\mu_{B^2,s^2h})}.
$$
Using (\ref{ineg84}), one obtains
$$
\begin{array}{llll}
 \displaystyle
 \int_0^1 |f(s)| ||  G(\widetilde{\pi_E}(X+s \cdot))-
\widetilde{G}(X+s \cdot) ||_{L^q(B^2,\mu_{B^2,h})} \ ds\\
\displaystyle \leq 
\int_0^1 |f(s)|  e^{\frac{|X|^2}{2h}} ||G||_{1,\varepsilon}
\sqrt{2hs^2} \left(\pi^{-\frac{1}{2}}
 \Gamma\left(\frac{q+s^{-2}+1}{2} \right) \right)^{\frac{1}{q+ s^{-2}}}
\sum_{0}^{\infty} \varepsilon_j( | u_j-\pi_E(u_j)| + | v_j-\pi_E(v_j)|)\ ds.
\end{array}
$$
For large values of $|z|$ and $|\arg(z)|<\pi$, one has (see \cite{MOS} 
p 12 for example)
$$
\Gamma(z)= z^{-\frac{1}{2}} e^{z(\ln(z)-1)} \sqrt{2\pi\ } (1+ O(z^{-1}) ).
$$
It follows that, for a constant $C$ which is independent of the parameters,
$$\left(\pi^{-\frac{1}{2}}
 \Gamma\left(\frac{q+s^{-2}+1}{2} \right) \right)^{\frac{1}{q+ s^{-2}}}
\leq C ( q+s^{-2}+1)^{\frac{1}{2}},
$$
which gives the estimate for the second term. \hfill $\square$

By differentiating under the integral sign one can obtain a weaker result,
 implying the existence of a stochastic extension, but not its precise form:
\begin{prop}
Let $G\in S_{m}(\bbb,\varepsilon)$ and let $f \in L^1([0,1])$, with values in
 $\R$. Let $X\in H^2$. For all  $Y\in H^2$, one sets
$$
T_XG(Y)= \int_0^1 f(s) G(X+sY)\ ds.
$$
Then $T_XG\in  S_{m}(\bbb,\varepsilon)$ and
$$
|| T_XG||_{m,\varepsilon}\leq  \int_0^1 |f(s)| ds \ || G||_{m,\varepsilon}.
$$
\end{prop}

\begin{coro}\label{moche}
Let $f$ be a continuous function on $[0,1]$, with values in  $\R$. Let
 $G$ be in $ S_1(\bbb,\varepsilon)$, with $\varepsilon$ summable.
Let $a_1,\dots, a_n$ and $X$ belong to $H^2$ and let
 $(p,h)$ be in $ [1,+\infty[\times \R^{+*}$. The function
$$
Y\in H^2  \mapsto
 \left( \prod_{i=1}^k<a_i,Y> \right) \ \int_0^1 f(s) G(X+sY)\ ds
$$
admits as a stochastic extension in $L^p(B^2,\mu_{B^2,h})$ the
 function
$$
Y \in B^2 \mapsto
 \left( \prod_{i=1}^k\ell_{a_i}(Y) \right) \
 \int_0^1 f(s) \widetilde{G}(X+sY)\ ds,
$$
where $ \widetilde{G}$ is the stochastic extension of $G$ valid for all
 $h' \leq h_0:=2h$ and all finite $p$.
Moreover, there exists a constant $K$ depending on  $p,k,h$ but not on the
 $a_i$,  $G,X,f,E$ or $\varepsilon$ such that, for all  $E\in \fff(H^2)$,
$$
\begin{array}{lll}
\displaystyle 
\left\Vert
 \displaystyle
 \left( \prod_{i=1}^k<a_i,\tilde{\pi}_E(Y)> \right)
 \ \int_0^1 f(s) G(X+s\tilde{\pi}_E(Y))\ ds
-  \left( \prod_{i=1}^k\ell_{a_i}(Y) \right) \
 \int_0^1 f(s) \widetilde{G}(X+sY)\ ds
 \right\Vert_{L^p(B^2,\mu_{B^2,h})}\\ \\
\leq \displaystyle 
K \int_0^1|f(s)|\ ds ||G||_{1,\varepsilon} A^{k-1}\times\\
\displaystyle \left(   \sum_{i=1}^k|\pi_E(a_i)-a_i| +
A\sqrt{\sum_{j=0}^{\infty}\varepsilon_j^2}
|\pi_E(X)-X| + A e^{|X|^2/2h}
 \sum_{j=0}^{\infty}\varepsilon_j\Big( |\pi_E(u_j)-u_j|+|\pi_E(v_j)-v_j|\Big)
 \right)
\end{array}
$$
where $A=\max_{1\leq i\leq k}(|a_i|)$.
\end{coro}
{\it  Proof.}
 % vérifiée le 13 fevrier 2016 mais jeté le papier et encore le 13 mars
One uses (\ref{HT-precise}) to establish that 
$$
\left\Vert
 \left( \prod_{i=1}^k<a_i,\tilde{\pi}_E(Y)> \right)
 \ \int_0^1 f(s) G(X+s\tilde{\pi}_E(Y))\ ds
-  \left( \prod_{i=1}^k\ell_{a_i}(Y) \right) \
 \int_0^1 f(s) \widetilde{G}(X+sY)\ ds
 \right\Vert_{L^p(B^2,\mu_{B^2,h})}$$
is smaller than
$$
\begin{array}{lll}
\displaystyle
\prod_{i=1}^k  ||\ell_{a_i}||\ \times \ 
 ||\int_0^1 f(s) G(X+s\tilde{\pi}_E(Y))\ ds- \int_0^1 f(s) \widetilde{G}(X+sY)\ ds ||\\
+ \displaystyle \sum_{i=1}^k\left(
 \prod_{j=1}^{i-1} ||\ell_{\pi_E(a_j)}||
\prod_{j=i+1}^{k} ||\ell_{a_j}||\right) \
 ||\ell_{\pi_E(a_i)}-\ell_{a_i}||\times 
 ||  \int_0^1 f(s) G(X+s\tilde{\pi}_E(Y))\ ds || , \\
\end{array}
$$
the norm in the second term being the  $L^{p(k+1)}(B^2,\mu_{B^2,h})$-norm.
But $||\ell_a|| =  K(p(k+1))h^{\frac{1}{2}}|a|$, according to (\ref{(44AJN)}) and
 (\ref{K(p)}). An upper bound is, consequently,
$$
\begin{array}{lll}
\displaystyle
(  K(p(k+1))h^{\frac{1}{2}})^k\left(  \prod_{i=1}^k  |a_j|\right)\ 
 ||\int_0^1 f(s) G(X+s\tilde{\pi}_E(Y))\ ds- \int_0^1 f(s)
 \widetilde{G}(X+sY)\ ds ||\\
+\displaystyle  (  K(p(k+1))h^{\frac{1}{2}})^k  \sum_{i=1}^k  |\pi_E(a_i)-a_i|
 \left(
 \prod_{1\leq j\leq k, j\neq i} |a_j|\right)\ 
\times 
 ||  \int_0^1 f(s)  G(X+s\tilde{\pi}_E(Y))\ ds ||  \\
\displaystyle
\end{array}
$$
One concludes by remarking that $|G|$ is smaller than  $||G||_{1,\varepsilon}$, 
(thanks to Proposition \ref{prol-int}) and that
the  $|a_i|$ are smaller than $A$. \hfill $ \square$

%%%%%%%%%%%%%%%%%%%%%%%%%%%%%%%%%%%%%%%%%%%%%%%%%%%%%%%%%%%%%%%%%%%%%%%%%
%%%%%%%%%%%%%%%%%%%%%%%%%%%%%%%%%%%%%%%%%%%%%%%%%%%%%%%%%%%%%%%%%%%%%%%%%
%%%%%%%%%%%%%%%%%%%%%%%%%%%%%%%%%%%%%%%%%%%%%%%%%%%%%%%%%%%%%%%%%%%%%%%%%
%%%%%%%%%%%%%%%%%%%%%%%%%%%%%%%%%%%%%%%%%%%%%%%%%%%%%%%%%%%%%%%%%%%%%%%%%

\section{Taylor expansions}\label{4-Taylor}  %4

\subsection{Differentiability of the symbols in $S_m(\bbb,\eps)$} %4.1
In this subsection, the sequence $\eps$ is supposed to be square summable
in most results.\\

The following straightforward lemma lists useful properties of the
 $S_m$ classes:
\begin{lemm}\label{Laplacien-H}
Let $F\in S_m(\bbb,\varepsilon)$, with $\varepsilon$ square summable.
\begin{itemize}
\item
If $m\geq 1$, for all  $i\in \Gamma$,  $\partial_{u_i} F$ and $\partial_{v_i} F$
belong to $ S_{m-1}(\bbb,\varepsilon)$ and
 $||\partial_{u_i} F||_{m-1,\eps}\leq \eps_i ||F||_{m,\eps}$, 
 $||\partial_{v_i} F||_{m-1,\eps} \leq \eps_i ||F||_{m,\eps}$.
More generally, if $m\geq k\geq 1$ and if $\alpha,\beta$ are two multi-indexes
 of depth $k$ (such that $\max_{j\in\Gamma}(\alpha_j,\beta_j)\leq k$), then 
$\partial^{\alpha}_u\partial^{\beta}_v F\in S_{m-k}(\bbb,\eps)$ and
$$
|| \partial^{\alpha}_u\partial^{\beta}_v F ||_{m-k,\eps} \leq ||F||_{m,\eps}
 \prod_{j\in\Gamma} \eps_j^{\alpha_j+\beta_j}.
$$
\item
If $m\geq 2$, one defines  $\Delta_{\bbb}$ by
$$
\Delta_{\bbb} F=\left( \sum_{j\in \Gamma} 
\left(\frac{\partial}{\partial u_j} \right)^2
+ \left(\frac{\partial}{\partial v_j} \right)^2\right) F.
$$
It is well defined and $\Delta_{\bbb} F\in S_{m-2}(\bbb,\varepsilon)$, with
$|| \Delta_{\bbb} F||_{m-2,\eps} \leq 2\sum_j\varepsilon_j^2|| F||_{m,\eps}$.
\item
If $G \in S_m(\bbb, \delta)$ with $\delta$ square summable too,
then $FG\in S_m(\bbb, \varepsilon + \delta)$ with
$|| FG||_{m,\varepsilon + \delta}\leq ||F||_{m,\varepsilon } ||G||_{m, \delta}$.
\end{itemize}
\end{lemm}
One can prove that, under certain conditions, the Laplace operator does
 not depend on the chosen 
basis (see Remark \ref{laplacien-independant} below).

\begin{prop}\label{Frechet-differentiable}
If $F\in S_m(\bbb,\varepsilon)$  with $m\geq 2$ and $\varepsilon$ 
 square summable, then $F$ is Fr\'echet differentiable on $H^2$ and
$$
DF(X)\cdot Y = \sum_{j\in \Gamma} <Y,u_j>\frac{\partial F}{\partial u_j}(X) + 
 <Y,v_j>\frac{\partial F}{\partial v_j}(X).
$$
Moreover, for all $X$ and $Y$ in $H^2$, 
$$
|F(X+Y)-F(X) - DF(X)\cdot Y|\leq ||F||_{m,\eps}
 \sum_{j\in \Gamma} \varepsilon_j^2(1+2\sqrt{2}) \ |Y|^2.
$$
\end{prop}
{\it Proof.}
Let $X,Y\in H^2$. Suppose that $\Gamma$ is enumerated. Let $P_N$ be the
 orthogonal projection  onto ${\rm Vect}(u_i,v_i, i\leq N)$ if $N\geq 0$,
 $P_{-1}=0$  and  $P_{N,\frac{1}{2}}$ the orthogonal projection  onto
 ${\rm Vect}(u_i,v_j,i\leq N+1, j\leq N)$. By approaching $P(X+Y)$ by
 $P(X+P_N(Y))$ one obtains
$$
\begin{array}{lll}
\displaystyle
F(X+P_N(Y)) - F(X)=\sum_{j=0}^N F(X+ P_j(Y)) - F(X+P_{j-1,\frac{1}{2}}(Y)) 
+F(X+P_{j-1,\frac{1}{2}}(Y)) -F(X+P_{j-1}(Y)).
\end{array}
$$
Taylor's formula gives, for example for the part of the $j$-th term concerned
 with $v_j$, 
$$
\begin{array}{lll}
\displaystyle
F(X+ <Y,v_j> v_j+ P_{j-1,\frac{1}{2}}(Y) ) - F(X + P_{j-1,\frac{1}{2}}(Y) )\\
= \displaystyle  <Y,v_j> \frac{\partial F}{\partial v_j}(X+ P_{j-1,\frac{1}{2}}(Y) ) \\
\displaystyle+ 
 <Y,v_j>^2 \int_0^1 (1-s) \frac{\partial^2 F}{\partial v_j^2}(X+ P_{j-1,\frac{1}{2}}(Y) +s <Y,v_j> v_j) \
 ds\\
 \displaystyle 
=  <Y,v_j> \frac{\partial F}{\partial v_j}(X )\\
\displaystyle +  
<Y,v_j> \left( \frac{\partial F}{\partial v_j}(X+ P_{j-1,\frac{1}{2}}(Y) )
-\frac{\partial F}{\partial v_j}(X )\right)\\
\displaystyle + 
 <Y,v_j>^2 \int_0^1 (1-s) \frac{\partial^2 F}{\partial v_j^2}(X+ P_{j-1,\frac{1}{2}}(Y) +s <Y,v_j> v_j) \
 ds.\\
\end{array}
$$
The first term gives the expression of the differential and it is the general 
term of a convergent series (apply Cauchy-Schwarz inequality). Since
 $\frac{\partial F}{\partial v_j}$ is in $S_{m-1}(\bbb,\varepsilon)$ with
$\displaystyle \left\Vert\frac{\partial F}{\partial v_j}\right\Vert_{m-1,\eps}
\leq
\eps_j ||F||_{m,\eps}$,
one can use (\ref{LIP}) to treat the second term. It then yields a convergent
series too, its sum being smaller than ${\rm Cste.}|Y|^2$. The integral term
can be estimated thanks to the estimates on the second derivatives and the sum
of the corresponding terms is of order $2$ in $|Y|$. Since $F$ and its
 derivatives are bounded by  $||F||_{m,\eps}$  and powers of $\eps$ 
independently on $X$ and $Y$, the rest can be bounded as is asserted in the
 theorem, with a constant $C$ independent of $X,Y$, $||F||_{m,\eps}$ and
 $\varepsilon$. One can take  $C=(1+2\sqrt{2})$.
\hfill $\square$\\

\noindent{\bf Remark.} 
Since there are infinitely many terms, we need a precise bound for the
rest in Taylor's formula, which explains the loss of one order of
 differentiability.\\
%%%%%%%%%%%%%%%%%%%%%%%%%%%%%%%%%%%%%%%%%%%%%%%%%%%%%%%%%%%%%%%%%%%%%%%%
%%%%%%%%%%%%%%%%%%%%%%%%%%%%%%%%%%%%%%%%%%%%%%%%%%%%%%%%%%%%%%%%%%%%%%%%

Deriving term by term and using the continuity of the extension operator
$\ppp$  (Corollary \ref{op-prol}) gives the following results:
\begin{prop}\label{Frechet-differentielle-Sm-1}
Let $F\in S_m(\bbb, \varepsilon)$  with $m\geq 2$,   $\varepsilon$ 
 square summable. Then, for all  $Y\in H^2$, 
$
X\mapsto DF(X)\cdot Y
$
is in $S_{m-1}(\bbb,\varepsilon),$ with
$|| X\mapsto DF(X)\cdot Y||_{m-1,\eps} \leq 
2||F||_{m,\eps}|Y|\sqrt{\sum_{j\in \Gamma}\varepsilon_j^2}$.
\end{prop}

%%%%%%%%%%%%%%%%%%%%%%%%%%%%%%%%%%%%%%%%%%%%%%%%%%%%%%%%%%%%%%%%%%%%%%%%
%%%%%%%%%%%%%%%%%%%%%%%%%%%%%%%%%%%%%%%%%%%%%%%%%%%%%%%%%%%%%%%%%%%%%%%%

\begin{coro}
Let  $F\in S_m(\bbb,\varepsilon)$ with $m\geq 2$ and $\varepsilon$ summable. 
The application $X\mapsto DF(X)\cdot Y $ from $H$ in  $\R$ admits a
 stochastic extension in $L^p(B^2,\mu_{B^2,t})$, which is the 
application
\begin{equation}\label{diff-somme}
\sum_{\Gamma} <Y,u_j> \ppp\left( \frac{\partial F}{\partial u_j}\right)
+  <Y,v_j> \ppp\left( \frac{\partial F}{\partial v_j}\right).
\end{equation}
\end{coro}
Here, the summability of $\eps$ is needed to ensure the existence of the 
stochastic extension.

%%%%%%%%%%%%%%%%%%%%%%%%%%%%%%%%%%%%%%%%%%%%%%%%%%%%%%%%%%%%%%%%%%%%%%%%
%%%%%%%%%%%%%%%%%%%%%%%%%%%%%%%%%%%%%%%%%%%%%%%%%%%%%%%%%%%%%%%%%%%%%%%%

\begin{defi}\label{DkF-forme}
Let $F\in S_m(\bbb,\varepsilon)$ with $\varepsilon$  square summable. For
 $k\in \{1,\dots, m\}$ and $X\in H^2$,  one defines a $k$-linear symmetric
continuous form  $\Phi_k(X)$ on $(H^2)^k$ setting:
$$
\begin{array}{lll}
\displaystyle
\forall (Y_1,\dots, Y_k) \in (H^2)^k,\\ \displaystyle \Phi_k(X)(Y_1,\dots,Y_k)=
& \displaystyle
\sum_{
\begin{array}{lll}(j_1,\dots, j_k)\in \Gamma^k,\\ (\delta_1,\dots,\delta_k)\in \{0,1\}^k 
\end{array}
}
\left( \prod_{s=1}^k <Y_s,w_{j_s}^{\delta_s}> \right)
\frac{\partial^k F}{\partial w_{j_1}^{\delta_1}\dots \partial w_{j_k}^{\delta_k}}(X),
\end{array}
$$
with $w_j^0=u_j, w_j^1=v_j$.
Moreover 
\begin{equation}\label{Cont-Phi-k}
\forall X,Y_1,\dots, Y_k \in ( H^2)^{k+1},\ 
| \Phi_k(X)(Y_1,\dots,Y_k)|\leq 2^k ||F||_{m,\eps} \prod_{s=1}^k |Y_s|
 \left( \sum_{\Gamma}\varepsilon_j^2\right)^{\frac{k}{2}}.
\end{equation}
\end{defi}

From now on, for the sake of brevity, we shall write
$J\in  \Gamma^k, \delta \in  \{0,1\}^k $ instead of
 $(j_1,\dots, j_k)\in \Gamma^k, (\delta_1,\dots,\delta_k)\in \{0,1\}^k $.

\begin{prop}
Let $F\in S_m(\bbb,M,\varepsilon)$ with $\varepsilon$  square summable. Then 
$F$ is  $C^{m-1} $ on  $H^2$ and, for all  
$k\in \{1,\dots, m-1\} $ and all $X\in H^2$, 
$$ D^kF(X)= \Phi_k(X).
$$
The inequality (\ref{Cont-Phi-k}) is satisfied. Finally, for $0\leq k\leq m-2$,
 one has
$$
||| D^kF(X+Z) -D^kF(X) - D^{k+1}F(X)(\cdot,Z)|||
\leq
2^k ||F||_{m,\eps}
\left(\sum_{\Gamma} \varepsilon_j^2\right)^{(k+2)/2} (1+2\sqrt{2}) |Z|^2,
$$
where the  norm is the norm of $k$-linear continuous applications on $H^2$. 
\end{prop}
{\it Proof.}
Propositions \ref{Frechet-differentiable} and \ref{Frechet-differentielle-Sm-1}
give the result for $m=2$, except for the fact that $F$ is $C^1$. This can
 be proved by applying (\ref{LIP}) to the partial derivatives of $F$. 
For a general $m$, one uses induction.
\hfill $\square$\\

This allows to state Taylor's formula to the order $k$ for  
 $F\in S_m(\bbb,\varepsilon)$, with  $\varepsilon$ square summable 
and $m\geq k+1$.
For $X, Y\in H^2$, 
\begin{equation}\label{Taylor}
\begin{array}{lll}
\displaystyle 
F(X+Y) & \displaystyle = F(X) +\sum_{i=1}^{k-1} \frac{1}{i!} D^iF(X)\cdot Y^i 
  +\displaystyle \int_0^1 \frac{(1-s)^{k-1}}{(k-1)!} D^kF(X+sY)\cdot Y^k\ ds
\\
 & \displaystyle = F(X) +\sum_{i=1}^{k-1} \frac{1}{i!} D^iF(X)\cdot Y^i \\
 & +\displaystyle \sum_{J\in \Gamma^k, \delta \in \{0,1\}^k }
\left( \prod_{r=1}^k <w_{j_r}^{\delta_r}, Y>\right)
\int_0^1 \frac{(1-s)^{k-1}}{(k-1)!} 
 \frac{\partial^k F}{\partial w_{j_1}^{\delta_1}\dots\partial w_{j_k}^{\delta_k}}
(X+sY) \ ds,
\end{array}
\end{equation}
exchanging the sums to get the last equality.\\

One part of the following subsection  \ref{subsection-Taylor-prolstoch} proves the existence of stochastic extensions for each of the terms
appearing here, the polynomial terms as well as the rest, under the 
assumption that $\eps$ is summable. Note that these extensions are 
 series indexed by   $\Gamma$.\\

One can finally state the following result, which allows us to construct 
another completion $B_A$ of $H$ in the case when  $\eps$ is
 summable.
\begin{prop}\label{BA-ancienneclasse}
Let  $\eps$ be a summable  sequence such that $\eps_j>0$ for all $j\in \Gamma$.
One defines a symmetric, definite positive and trace class operator $A$
by setting 
$$
\forall X\in B^2, \ AX= \sum_{j\in \Gamma} 
\eps_j<X, u_j>u_j + \eps_j<X, v_j> v_j.
$$
Set $||X||_A= <AX, X>^{1/2}$. Then $|| \ ||_A$ is a measurable norm on $H$,
 in the sense of  \cite{K} (Def.4.4) or \cite{G-1}. 
One denotes by $B_A$ the completion of $H$ for this norm.
\\
If $F\in S_m(\bbb,\eps)$ for $m\geq 2$, then  $F$ is uniformly
continuous on $H^2$ with respect to the norm  $||\ ||_A$.
 The function $F$ admits a uniformly continuous  extension $F_A$ on $B_A$ and 
the stochastic extension $\tilde{F}$ of $F$  given by
 Proposition \ref{prop84-modifiee}
is  equal to $F_A$ $\mu_{B,h}$- a.e.
\end{prop}

 {\it Proof.}
It follows from Theorem 3 in  \cite{G-1} that $|| \ ||_A$ is a measurable norm,
since $A$ is injective. Since  $m\geq 2$, $F$ is $C^1$ on $H$. Taylor's formula
with an integral rest and Definition \ref{DkF-forme} 
allow us to write the inequality 
$$
\begin{array}{lll}
|F(X)-F(Y)| & \displaystyle \leq \int_0^1 \sum_{j\in\{0,1\},\delta\in \Gamma}
\left|\frac{\partial F}{\partial w_j^{\delta}}(X+t(Y-X)) <Y-X, w_j^{\delta}> 
\right| \ dt \\
& \leq \sum_{j\in \Gamma} ||F||_{m,\eps}\eps_j^{1/2} \ \eps_j^{1/2} ( |<Y-X, u_j>| + 
|<Y-X, v_j>|)   \\
& \displaystyle \leq || F||_{m,\eps} \sqrt{2}
 \left( \sum \eps_j \right)^{1/2} ||X-Y||_A,\\
\end{array}
$$
thanks to Cauchy-Schwarz inequality. This proves that $F$ is
 uniformly continuous on $H^2$ and therefore admits an
 extension $F_A$, which is uniformly continuous on $B_A$.
According to Theorem 6.3 (Chap 1 \cite {K}), $F_A$ and $\tilde{F}$
coincide almost everywhere.
\hfill $\square$

\begin{rem} \label{laplacien-independant}
If $F\in S_m(\bbb,\eps)$ with $m\geq 3$ and $\eps$ summable, one can define
$\Delta F$ more intrinsically. Indeed, one can state an  inequality 
more precise than (\ref{Cont-Phi-k}).  For  $k\leq 3$ one gets
$$
\forall X,Y_1,\dots, Y_k \in ( H^2)^{k+1},\ 
| \Phi_k(X)(Y_1,\dots,Y_k)|\leq 2^k ||F||_{m,\eps}
 (\sum_{\Gamma}\varepsilon_j)^{k/2}
\prod_{s=1}^k <AY_s,Y_s>^{1/2},
$$
reasoning as in the proof of  Proposition \ref{BA-ancienneclasse}. 
The function $F$ is $C^2 $ since  $m\geq 3$ and the
 inequality, for $k=2$, ensures the existence of a self adjoint, trace class
 operator $M_x$ satisfying
$$
\forall U,V,X \in H^2, d^2F(X)\cdot (U,V) = <M_XU,V>.
$$
One then sets $\Delta F(X)= {\rm Tr}(M_X)$ and the expression as a sum of
 partial derivatives does not depend on the chosen orthonormal basis.\\
One can remark, too, that if $\eps$ is summable, if $F$ belongs to
 $ S_m(\bbb,\eps)$ for all  $m$ and if there exists a constant $M$ such that
  $||F||_{m,\eps}\leq M$ for all  $m$, then
$F\in S(Q_B)$ with $B$ defined by  $B=4 (\sum_{\Gamma}\varepsilon_j) A$,
$A$ being as in Proposition \ref{BA-ancienneclasse}.
% for the snark was a boojum, you see
\end{rem}

%%%%%%%%%%%%%%%%%%%%%%%%%%%%%%%%%%%%%%%%%%%%%%%%%%%%%%%%%%%%%%%%%%%%%%%%
%%%%%%%%%%%%%%%%%%%%%%%%%%%%%%%%%%%%%%%%%%%%%%%%%%%%%%%%%%%%%%%%%%%%%%%%
%%%%%%%%%%%%%%%%%%%%%%%%%%%%%%%%%%%%%%%%%%%%%%%%%%%%%%%%%%%%%%%%%%%%%%%%
%%%%%%%%%%%%%%%%%%%%%%%%%%%%%%%%%%%%%%%%%%%%%%%%%%%%%%%%%%%%%%%%%%%%%%%%
%%%%%%%%%%%%%%%%%%%%%%%%%%%%%%%%%%%%%%%%%%%%%%%%%%%%%%%%%%%%%%%%%%%%%%%%
\subsection{Taylor's formula and stochastic expansions}
\label{subsection-Taylor-prolstoch} % 4.2

Contrary to the preceding subsection, where sums like
 $\sum_{\Gamma} \eps_j< u_j,x>$ have been treated by Cauchy-Schwarz inequality,
we must suppose here that the sequence $\eps$ is summable. The corresponding
sums have the form  $\sum_{\Gamma} \eps_j\ \ell_{u_j}$ and, since 
the functions $ \ell_{u_j}$ have a $L^p$ norm independent of $j$, 
 Cauchy-Schwarz inequality cannot be applied. 

\begin{lemm}
Let  $\varepsilon$ be a summable sequence. Let  $F\in S_m(\bbb,\varepsilon)$
 with $m\geq 2$ and let  $X\in H^2$. For all  $k\leq m$ and all $h>0$,
 $p\in [1,+\infty[$, the application $Y\mapsto \Phi_k(X)\cdot Y^k$ from
 Definition  \ref{DkF-forme}  admits, as a stochastic expansion in
 $L^p(B^2,\mu_{B^2,h})$, the application
$Y\mapsto \widetilde{ \Phi_k(X)}\cdot Y^k$ defined on $B^2$ by
$$
\forall Y \in B^2,\qquad 
\widetilde{ \Phi_k(X)}\cdot Y^k \ =\ 
\sum_{J\in \Gamma^k, \delta \in \{0,1\}^k }
\left( \prod_{s=1}^k \ell_{w_{j_s}^{\delta_s}}( Y)\right)
 \frac{\partial^k F}{\partial w_{j_1}^{\delta_1}\dots\partial w_{j_k}^{\delta_k}}(X),
$$
with $w_j^0=u_j, w_j^1=v_j$.
\end{lemm}
{\it Proof.}
Let $E\in \fff(H)$. To verify that  $\widetilde{ \Phi_k(X)}\cdot Y^k$
and $Y\mapsto \Phi_k(X)\cdot(\tilde{\pi}_E(Y)^k)$  really belong to
 $L^p(B^2,\mu_{B^2,h})$, one has to find an upper bound for each term
$$
\left\Vert  \prod_{s=1}^k \ell_{a_s}( Y)  \right\Vert_{L^p(B^2,\mu_{B^2,h})}
\left| \frac{\partial^k F}{\partial w_{j_1}^{\delta_1}\dots
\partial w_{j_k}^{\delta_k}}(X) \right|
$$
of the sum, with $a_s= w_{j_s}^{\delta_s} $ ou $\pi_E(w_{j_s}^{\delta_s})$.
One then proves, using Proposition \ref{polyn-prodscal}, that
$$
\left\Vert \prod_{s=1}^k <\tilde{\pi}_E(Y), w_{j_s}^{\delta_s}> - 
\prod_{s=1}^k \ell_{w_{j_s}^{\delta_s}} \ \right\Vert_{L^p(B^2,\mu_{B^2,h})}
\ \leq \ 
(K(pk)h^{\frac{1}{2}} )^{k} \sum_{s=1}^k 
|\pi_E(w_{j_s}^{\delta_s}) - w_{j_s}^{\delta_s}|,
$$
since the $w_{j_s}^{\delta_s}$  and their projections have norms smaller than $1$
Therefore
$$
\begin{array}{llll}
\displaystyle 
\Vert  \Phi_k(X)\cdot(\tilde{\pi}_E(Y)^k)-\widetilde{ \Phi_k(X)}(Y,\dots,Y)
\Vert_{L^p(B^2,\mu_{B^2,h})}\\
\leq \displaystyle ||F||_{m,\eps} (K(pk)h^{\frac{1}{2}})^{k}
\sum_{J\in \Gamma^k, \delta \in \{0,1\}^k }
 \prod_{s=1}^k \varepsilon_{j_s} \sum_{s=1}^k
|\pi_E(w_{j_s}^{\delta_s}) - w_{j_s}^{\delta_s}|.
\end{array}
$$
One then replaces $E$ by $E_n$, where $(E_n)$ is an increasing sequence  of
 $\fff(H^2)$ whose union is dense in $H^2$. Since the terms 
$|\pi_{E_n}(w_{j_s}^{\delta_s}) - w_{j_s}^{\delta_s}|$ converge to $0$ 
and are smaller than $2$, the difference converges to $0$ thanks to
 the dominated convergence Theorem. \hfill $\square$\\

\begin{prop}\label{prop-Tayloretendue}
Let $\varepsilon$  be summable and let $F\in S_m(\bbb,\varepsilon)$ with
  $m\geq 2$. Let $X\in H^2$. For all  $k\leq m-1$,
 all $h>0$ and $p\in [1,+\infty[$, one can write, in $L^p(B^2,\mu_{B^2,h})$:
\begin{equation}\label{Tayloretendue}
\begin{array}{lll}
\displaystyle 
\widetilde{F}(X+Y) & \displaystyle = F(X) 
+  \sum_{i=1}^{k-1} \frac{1}{i!}
\sum_{J\in \Gamma^i, \delta \in \{0,1\}^i }
\left( \prod_{r=1}^i \ell_{w_{j_r}^{\delta_r}}(Y)\right)
 \frac{\partial^i F}{\partial w_{j_1}^{\delta_1}\dots\partial w_{j_i}^{\delta_i}}(X)
\\
 & +\displaystyle 
\sum_{J\in \Gamma^k, \delta \in \{0,1\}^k }
\left( \prod_{r=1}^k \ell_{w_{j_r}^{\delta_r}}(Y) \right)
\int_0^1 \frac{(1-s)^{k-1}}{(k-1)!} 
 \ppp\left(\frac{\partial^k F}{\partial w_{j_1}^{\delta_1}\dots\partial w_{j_k}^{\delta_k}}\right)
(X+sY) \ ds.
\end{array}
\end{equation}
\end{prop}

{\it Proof.}
Let us denote by $\widetilde{\Phi}_i(X)\cdot Y^i$ the $i$-th term of the sum
 and by  $R_k(X)$ the last one, corresponding to the rest.
We have just seen that the polynomial part of the development in
 (\ref{Taylor}) has a stochastic extension in $L^p(B^2,\mu_{B^2,h})$.
The rest is the sum indexed by
 $J=(j_1,\dots, j_k)\in \Gamma^k,
\delta= (\delta_1,\dots,\delta_k)\in \{0,1\}^k $. 
One applies the Corollary \ref{moche}, replacing, in the upper bound, 
$\int_0^1 \frac{(1-s)^{k-1}}{(k-1)!} \ ds$ by $(k!)^{-1}$,
 $||G||_{1,\varepsilon}$ by $||F||_{m,\varepsilon}\prod_{1}^k \varepsilon_{j_i}$
 and $A=\max(| w_{j_i}^{\delta_i}|)$, by $1$. One finds
$$
\begin{array}{lll}
\displaystyle \sum_{J\in \Gamma^k, \delta \in \{0,1\}^k } &\displaystyle
\Big\Vert \int_0^1  \frac{(1-s)^{k-1}}{(k-1)!} 
\Big(
 \frac{\partial^k F}{\partial w_{j_1}^{\delta_1}\dots\partial w_{j_k}^{\delta_k}}
(X+s\tilde{\pi}_E(Y)) \prod_1^k <\tilde{\pi}_E(Y),w_{j_i}^{\delta_i} >
\\
&\displaystyle 
-(\ppp\frac{\partial^k F}{\partial w_{j_1}^{\delta_1}\dots\partial w_{j_k}^{\delta_k}})
(X+sY)  \prod_1^k \ell_{w_{j_i}^{\delta_i}}(Y)  \Big)\ ds
 \Big\Vert_{L^p(B^2,\mu_{B^2,h})} \\
&\leq \displaystyle 
\frac{1}{k!} K ||F||_{m,\varepsilon}  \sum_{J\in \Gamma^k, \delta \in \{0,1\}^k }
  (\varepsilon_{j_1} \dots \varepsilon_{j_k}) \sum_{i=1}^k
 |\pi_E(w_{j_i}^{\delta_i} )-w_{j_i}^{\delta_i} |           \\
&+ \displaystyle 
\frac{1}{k!} K ||F||_{m,\varepsilon}  \sum_{J\in \Gamma^k, \delta \in \{0,1\}^k } 
 (\varepsilon_{j_1} \dots \varepsilon_{j_k})       \quad  \times   \\
&\displaystyle \left( \sqrt{\sum\varepsilon_j^2} |\pi_E(X)-X|
+e^{|X|^2/2h}
\sum_0^{\infty} \varepsilon_j(|\pi_E(u_j)-u_j| +|\pi_E(v_j)-v_j|) \right).
\end{array}
$$
If one replaces $E$ by  $E_n$ from an increasing sequence of $\fff(H^2)$
whose union is dense in $H^2$, this  converges to  $0$ when $n$ converges to
 infinity.\hfill $\square$\\

%%%%%%%%%%%%%%%%%%%%%%%%%%%%%%%%%%%%%%%%%%%%%%%%%%%%%%%%%%%%%%%%%%%
%%%%%%%%%%%%%%%%%%%%%%%%%%%%%%%%%%%%%%%%%%%%%%%%%%%%%%%%%%%%%%%%%%%
%%%%%%%%%%%%%%%%%%%%%%%%%%%%%%%%%%%%%%%%%%%%%%%%%%%%%%%%%%%%%%%%%%%

With each term of the extended Taylor expansion (\ref{Tayloretendue}), 
one can associate a quadratic form (see \cite{AJN}, Definition 1.2)
thanks to the following result:
\begin{prop}\label{4.11}
Let $F\in S_m(\bbb,\varepsilon)$ with $\varepsilon$ summable and $m\geq k+1$,
 where $k$ is the  order of differentiation. Each of the terms of 
 (\ref{Tayloretendue}) has a  $N_s$ norm (cf. (\ref{(1.12AJN)})), for a
 well-chosen $s$. Precisely
$$
N_i(\frac{1}{i!} \widetilde{\Phi}_i(X)\cdot Y^i) \leq \frac{1}{i!}
||F||_{m,\varepsilon}
 \left(2 \max(1,\sqrt{\frac{h}{2}})\sum_{\Gamma}\varepsilon_{j}\right)^i
\int_{\R} (1+|v|)^i d\mu_{\R,1}(v) .
$$
and
$$
N_k( R_k(X)  ) \leq \frac{1}{k!}
||F||_{m,\varepsilon}
 \left(2 \max(1,\sqrt{\frac{h}{2}})\sum_{\Gamma}\varepsilon_{j}\right)^k
\int_{\R} (1+|v|)^k d\mu_{\R,1}(v) .
$$
\end{prop}
{\it Proof.}
One uses the computations of Proposition \ref{NormeNm-fonctionsell}. Then
$$
|| \prod_{r=1}^i \ell_{w_{j_r}^{\delta_r}}(\cdot +Y)||_{L^1(B^2,\mu_{B^2,\frac{h}{2}})}
\leq (1+|Y|)^i \max(1,\sqrt{\frac{h}{2}})^i 
\int_{\R} (1+|v|)^i d\mu_{\R,1}(v).
$$
Hence
$$
\begin{array}{lll}
\displaystyle
\left\Vert \sum_{J\in \Gamma^i, \delta \in \{0,1\}^i }
\prod_{r=1}^i \ell_{w_{j_r}^{\delta_r}}(\cdot +Y) 
\frac{\partial^i F}{\partial w_{j_1}^{\delta_1}\dots
\partial w_{j_i}^{\delta_i}}(X)\right\Vert_{L^1(B^2,\mu_{B^2,\frac{h}{2}})}\\
\displaystyle
\leq \sum_{J\in \Gamma^i, \delta \in \{0,1\}^i }
 ||F||_{m,\varepsilon}\varepsilon_{j_1}\dots
\varepsilon_{j_i}\  (1+|Y|)^i \max(1,\sqrt{\frac{h}{2}})^i 
\int_{\R} (1+|v|)^i d\mu_{\R,1}(v)\\
\displaystyle
\leq 
 ||F||_{m,\varepsilon}
 \left(2 \max(1,\sqrt{\frac{h}{2}})\sum_{\Gamma}\varepsilon_{j}\right)^i
\int_{\R} (1+|v|)^i d\mu_{\R,1}(v)  (1+|Y|)^i .\\
\end{array}
$$
It follows that
$$
N_i(\frac{1}{i!} \widetilde{\Phi}_i(X)\cdot Y^i) \leq \frac{1}{i!}
||F||_{m,\varepsilon}
 \left(2 \max(1,\sqrt{\frac{h}{2}})\sum_{\Gamma}\varepsilon_{j}\right)^i
\int_{\R} (1+|v|)^i d\mu_{\R,1}(v) .
$$
We treat the rest in the same way: the sum indexed by
 $(j_1,\dots, j_k)\in \Gamma^k, (\delta_1,\dots,\delta_k)\in \{0,1\}^k $
contains a product of $k$ terms $\ell$ and the integral, which is
bounded by
$\frac{1}{k!} ||F||_{m,\varepsilon}\varepsilon_{j_1}\dots \varepsilon_{j_k} $.
Therefore the rest has a  $N_k$ norm bounded like the polynomial terms.
 \hfill $\square$\\

%%%%%%%%%%%%%%%%%%%%%%%%%%%%%%%%%%%%%%%%%%%%%%%%%%%%%%%%%%%%%%%%%%%%%
%%%%%%%%%%%%%%%%%%%%%%%%%%%%%%%%%%%%%%%%%%%%%%%%%%%%%%%%%%%%%%%%%%%%%
%%%%%%%%%%%%%%%%%%%%%%%%%%%%%%%%%%%%%%%%%%%%%%%%%%%%%%%%%%%%%%%%%%%%%

\section{The heat operator on $H$}\label{5-chaleur}%5

\subsection{Definition}
The heat operator defined below associates a function defined on a (real,
 separable, infinite dimensional) Hilbert space, with a function defined on
 the same Hilbert space. We aim at extending the notion of heat operator,
 which is classical in the finite dimesional setting. The results proved here 
are different from the results obtained by  (\cite{K}, \cite{G-4}), inasmuch as
they are concerned with functions initially defined on 
$H$ (or $H^2$) and not on $B$.
\begin{defi}
Let $F$ be a function defined on $H$, admitting a stochastic extension in
 $L^p(B, \mu_{B,t})$ for a given $p\in [1,+\infty[$. One defines 
 $H_tF$ on $H$ by
\begin{equation}\label{def-op-chaleur}
 (H_tF)(X)= 
\int_{B} \widetilde{F}(X+Y) \ d\mu_{B,t}(Y)=
\int_{B} \widetilde{F}(Y) e^{- \frac{|X|^2}{2t}} e^{\ell_X/t} \ d\mu_{B,t}(Y),
\end{equation}
the second identity coming from (\ref{(46AJN)}).\\
If  $F$ is defined on the product $H^2$, one replaces $H$ by $H^2$ and
 $B$ by $B^2$.
\end{defi}

\begin{rem}\label{indep}
This definition does not depend on the stochastic extension chosen, nor on
 the measurable norm and on the completion of $H$ associated with it. Indeed,
 the fact that a sequence $F\circ \tilde{\pi}_{E_n}$ is a Cauchy sequence in
 $L^p(B,\mu_{B,h})$ is expressed by integrals on finite dimensional subspaces
 of $H$ (using (\ref{7AJN})) and not at all by integrals on $B$. Likewise, the 
integral of (\ref{def-op-chaleur}) does not depend on the integration space 
$B$, since it is a limit of integrals on finite dimensional spaces of $H$.
\end{rem}

\begin{prop}\label{semi-groupe!}
Let $F$ belong to a class $S(Q_A)$ of Definition \ref{Classe-DD-def} or to 
a class $S_m(\bbb,\eps)$, with $\eps$ summable, of Definition
\ref{ancienne-classe} .The semigroup property is verified:
for all positive $s,t$ and all $X$ in the Hilbert space, 
$$
H_t(H_sF)(X)= H_{t+s}F(X).
$$
 Moreover, one has (according to whether $F\in S(Q_A)$ or $S_m(\bbb,\eps)$,
\begin{equation}\label{ineg-chaleur-simple}
\forall X\in H^2,\
\left|(H_t F)(X) \right|\leq ||F||_{m,\eps} \ {\rm or}\ 
\forall X\in H,\
\left|(H_t F)(X) \right|\leq ||F||_{Q_A}.
\end{equation}
\end{prop}
{\it Proof.} 
We give the proof in the case when $F\in S(Q_A)$. Let  $B_A$ be the completion 
of $H$ with respect to the  measurable norm $||\ ||_A$ given by Proposition
 \ref{BA-classeDD}. The function $F$ is uniformly continuous on $H$ and extends 
continuously  as a function denoted by $F_A$, uniformly continuous and
bounded on $B_A$. By Theorem 6.3 (Chap 1) of \cite{K}, every stochastic
 extension of $F$ in $L^p(B_A,\mu_{B_A,h})$ coincides with $F_A$ $\mu_{B_A,h}$-a.e.
 One can thus, considering that the heat operator is being defined by 
integrating on $B_A$, write that
$$
\forall X\in H,\ H_tF(X)= \int_{B_A} F_A(X+Y) \ d\mu_{B_A,t}(Y).
$$
This formula allows  us to define a function, denoted by  $H_tF_A$, on $B_A$.
Since $F_A$ is uniformly continuous and bounded on $B_A$, $H_tF_A$ is uniformly
 continuous and bounded on $B_A$ too, by \cite{K} (Theorem 4.1 Chap 3). Then 
  $H_tF_A$ is the stochastic extension of its restriction to $H$, $H_tF$ and
$$
\forall X \in H,\ H_s(H_tF)(X)= \int_{B_A}H_tF_A(X+Y) 
d\mu_{B_A,s}(Y)= H_{t+s}F_A(X)= H_{t+s}F(X).
$$

For $F\in S_m(\bbb,\eps)$ with $\eps$ summable we can reproduce the 
same demonstration, with $H^2$ and  $||\ ||_A$, $B_A$ from 
Proposition \ref{BA-ancienneclasse}.\\
The inequalities   (\ref{ineg-chaleur-simple}) come from the fact that $F_A$ 
is bounded on $B_A$ like  $F$ on $H$.
\hfill $\square$

%%%%%%%%%%%%%%%%%%%%%%%%%%%%%%%%%%%%%%%%%%%%%%%%%%%%%%%%%%%%%%%%%%%%%%%%%%%%%%
%%%%%%%%%%%%%%%%%%%%%%%%%%%%%%%%%%%%%%%%%%%%%%%%%%%%%%%%%%%%%%%%%%%%%%%%%%%%%%
%%%%%%%%%%%%%%%%%%%%%%%%%%%%%%%%%%%%%%%%%%%%%%%%%%%%%%%%%%%%%%%%%%%%%%%%%%%%%%

\subsection{The heat operator  in the classes $S_m(\bbb,\eps)$}
\label{chaleur-anciennesclasses}

\begin{prop}\label{op-chaleur-H}
Let $F\in  S_m(\bbb,\varepsilon)$ with $m\geq 2$, $\varepsilon $
 summable. If $\alpha,\beta$ are depth $1$ multiindices (such that
 $\max(\alpha_j,\beta_j)\leq 1$), then
$$
\partial^{\alpha}_u \partial^{\beta}_v (H_t  F)(X) =
H_t(\partial^{\alpha}_u \partial^{\beta}_v)(X). 
$$
Moreover, for $m\geq 1$, $H_t F \in S_{m-1}(\bbb,\varepsilon)$, with
$\displaystyle || H_t  F ||_{m-1,\eps}\leq || F ||_{m,\eps} $.
The operator $H_t$ is continuous from  $S_m(\bbb,\varepsilon)$, in
 $S_{m-1}(\bbb,\varepsilon)$.
\end{prop}

{\it Proof.}
 If $m=1$, the continuity of $H_t$ from $S_1(\bbb,\eps)$ in  $S_0(\bbb,\eps)$
comes from the inequalities  (\ref{ineg-chaleur-simple}). Now suppose that
 $m\geq 2$ and prove (first) that
$$
\frac{\partial}{\partial w}(H_t  F)(X) =
H_t \left( \frac{\partial}{\partial w} F\right)(X) 
$$
 with $w=u_i$ or $v_i$ and $X\in H^2$. By Taylor's formula
\begin{equation}\label{Taylor2}
F(X+rw)-F(X)= r\frac{\partial F}{\partial w}(X)+ r^2 \int_0^1 (1-s) 
\frac{\partial^2 F}{\partial w^2}(X+rsw)\ ds.
\end{equation}
According to Proposition \ref{prop84-modifiee} and its corollary, $F$ and 
$\frac{\partial F}{\partial w}$ together with their translated of a vector 
$Y\in H^2$ admit stochastic extensions in $L^p(B^2,\mu_{B^2,t})$
and $\widetilde{\tau_Y F}= \tau_Y \widetilde{F}$.
According to (\ref{Taylor2}), for all  $r\in\R^*$, the function
$ G_r : X\mapsto  \int_0^1 (1-s) 
\frac{\partial^2 F}{\partial w^2}(X+rsw)\ ds $ 
 admits a stochastic extension in $L^p(B^2,\mu_{B^2,t})$, denoted by
 $\widetilde{G_r}$.\\
For all  $r$, $|G_r| \leq \frac{1}{2} ||F||_{m,\eps}\sup(\varepsilon_i)^2$.
 Hence, so does $\widetilde{G_r}$   $\mu_{B^2, t} -$ a.s.\\
Applying (\ref{Taylor2}) in the point  $\widetilde{\pi}_{E_j}(X)$ with 
$ X\in B^2$ and taking a limit in  $L^p(B^2,\mu_{B^2,t})$, one obtains 
\begin{equation}\label{Taylor-prolongement}
\tau_{rw}\widetilde{F}-\widetilde{F}= r{\cal P}(\frac{\partial F}{\partial w})
+ r^2\widetilde{G_r},\quad {\rm in }\  L^p(B^2,\mu_{B^2,t})
\end{equation}
One deduces that, for all  $X$ of $H^2$,
$$
\frac{({H_t}{F})(X+rw) 
-({H_t}{F})(X)}{r}=
\left({H_t}\frac{\partial F}{\partial w}\right)(X)
+ r({H_t}{G_r})(X),
$$
and that
$$
\left| \frac{(H_tF)(X+rw) - (H_tF)(X) }{r}
-(H_t\frac{\partial F}{\partial w})(X) \right|
\leq |r|\int_{B^2}|\widetilde{G_r}|(X+Y) d\mu_{B^2,t}(Y).
$$
The bound on $\widetilde{G_r}$ shows that
$$
\lim_{r\rightarrow 0}  \frac{(H_tF)(X+rw) - (H_tF)(X) }{r}=
 (H_t\frac{\partial F}{\partial w})(X),
$$
which means that $H_tF$ admits order $1$ partial derivatives in the (canonical)
 directions $u_i, v_i$ .\\
Let $\alpha,\beta$  be two depth $1$ multiindices.
Let $w=u_i$ (or $v_i$) be a coordinate, with respect to which one has not yet 
differentiated  (that is, such that  $\alpha_i =0$ or $\beta_i=0$). Applying
the preceding reasoning to $\partial^{\alpha}_u \partial^{\beta}_v F$, we get that
$$
\frac{\partial}{\partial w} H_t( \partial^{\alpha}_u \partial^{\beta}_v   F)(X)=
 H_t(\frac{\partial}{\partial w} \partial^{\alpha}_u \partial^{\beta}_v   F)(X)
$$
and an induction on $|\alpha|+|\beta|$ allows us to exchange  $H_t$ and
 differentiations. By (\ref{ineg-chaleur-simple}), one gets that 
$$
|\partial^{\alpha}_u \partial^{\beta}_v  H_t(   F)(X)|= 
| H_t( \partial^{\alpha}_u \partial^{\beta}_v   F)(X)|\leq 
||\partial^{\alpha}_u \partial^{\beta}_v   F||_{m-1,\eps}\leq
 \eps^{\alpha+\beta} 
||F||_{m,\eps}.
$$
If $m=2$, the proposition is proved. Otherwise one completes the proof by
 induction.\hfill $\square$ 
%%%%%%%%%%%%%%%%%%%%%%%%%%%%%%%%%%%%%%%%%%%%%%%%%%%%%%%%%%%%%%%%%%%%%%%
%%%%%%%%%%%%%%%%%%%%%%%%%%%%%%%%%%%%%%%%%%%%%%%%%%%%%%%%%%%%%%%%%%%%%%%
%%%%%%%%%%%%%%%%%%%%%%%%%%%%%%%%%%%%%%%%%%%%%%%%%%%%%%%%%%%%%%%%%%%%%%%

The Heat operator commutes with the Laplace operator:
\begin{prop}\label{anc-laplaciencomm}
Let $\varepsilon$ be summable. The operator $\Delta_{\bbb}$ is continuous from
$S_m(\bbb,\varepsilon)$  to $S_{m-2}(\bbb,\varepsilon)$, for $m\geq 2$.
Moreover, for   $m\geq 3$, 
$$
\forall F \in S_m(\bbb,\varepsilon),\quad
\Delta_{\bbb} H_t F = H_t \Delta_{\bbb} F \in S_{m-3}(\varepsilon).
$$
\end{prop}
{\it Proof.}
One deduces from Lemma  \ref{Laplacien-H}  that
$$
|| \Delta_{\bbb} F||_{m-2,\varepsilon}\leq  2 \sum_{j\in \Gamma} \varepsilon_j^2 
|| F||_{m,\varepsilon},
$$
which proves the continuity of $\Delta_{\bbb}$.\\
One still supposes  $\Gamma$ enumerated. For  $n\in \N$, set
 $\Delta_n= \sum_{j\leq n} 
\frac{\partial^2}{\partial u_j^2} +\frac{\partial^2}{\partial v_j^2}$. 
One can see that $\Delta_n F$ converges to  $\Delta_{\bbb} F$ in 
$S_{m-2}(\bbb,\varepsilon)$. Moreover, one can exchange $H_t$ and the
 differentiations with respect to $u_j, v_j$. This fact, and the continuity 
of the operators, allow us to write
$$
 H_t \Delta_{\bbb} F =
 H_t \lim_{n\rightarrow \infty} \Delta_n F =
 \lim_{n\rightarrow \infty} H_t \Delta_n F =
 \lim_{n\rightarrow \infty}\Delta_n H_t  F = \Delta_{\bbb} H_t  F,
$$
which completes the proof.
\hfill $\square$\\

Let us state a result about commutators. For $Z\in H^2$ and  $F$ a function
 defined on $H^2$, denote by $M_ZF$ the function defined by 
$
(M_ZF)(X)= <Z,X> F(X).
$

\begin{prop}
Let  $F\in S_m(\bbb,\varepsilon)$ with $m\geq 2$ and $\varepsilon $
square summable. For all  $i\in \N$, one has
$$
\frac{1}{t}\left(H_tM_{u_i}-M_{u_i} H_t\right )F = 
 H_t\frac{\partial F}{\partial u_i}
$$
and then
$$
\frac{1}{t}\left[H_t,M_{u_i}\right] = 
 H_t\frac{\partial }{\partial u_i}
$$
on $S_m(\bbb,\varepsilon)$. The same property holds with  $v_i$.
\end{prop}{\it Proof.}
Notice that $M_Z F$ admits $\ell_Z \widetilde{F}$ as a  stochastic extension
in $L^p(B^2, \mu_{B^2,t})$ for all  $p\in[1,+\infty[$, by Corollary \ref{corHT}.
According to Theorem 6.2 (chap. 2, par. 6) of \cite{K}, for all  $X\in H^2$, 
$$
\frac{\partial H_tF }{\partial u_i}(X)=
\frac{1}{t}
\int_{B^2} \widetilde{F}(X+Y) \ell_{u_i}(Y) \ d \mu_{B^2,t}(Y).
$$
But $\ell_{u_i}(Y)=\ell_{u_i}(Y+X)- <u_i,X>$, since $X\in H^2$. Then
$$
\frac{\partial H_tF }{\partial u_i}(X)=
\frac{1}{t}
\int_{B^2} \widetilde{F}(X+Y) \ell_{u_i}(Y+X) \ d \mu_{B^2,t}(Y)
- <u_i,X> \frac{1}{t}
\int_{B^2} \widetilde{F}(X+Y)  \ d \mu_{B^2,t}(Y).
$$
This is the desired result. \hfill $\square$\\

We shall use the Taylor expansions and their stochastic extensions to prove
a preliminary result before stating the main result of this subsection,
Theorem \ref{glavni-Sm}.
\begin{prop}\label{auxiliaire}
\begin{enumerate}
\item
Let $m\geq 3$. There exists  $C_m\in \R^+$ such that, for all  
 $F\in S_m(\bbb,\varepsilon)$,
\begin{equation}\label{cont-0}
\left\Vert H_tF -F\right\Vert_{m-3,\eps} \leq C_m ||F||_{m,\eps} t.
\end{equation}
For all  $s>0$, for $m\geq 5$, one has
\begin{equation}\label{contpartout}
\left\Vert H_{t+s}F -H_sF\right\Vert_{m-4,\eps} \leq C_m ||F||_{m,\eps} t.
\end{equation}
\item
Let $m\geq 4$. There exists $C_m\in \R^+$  such that, for all 
  $F\in S_m(\bbb,\varepsilon)$,
\begin{equation}\label{un}
\left\Vert \frac{H_tF -F}{t} -
 \frac{1}{2} \Delta F\right\Vert_{m-4,\eps} \leq C_m ||F||_{m,\eps} t^{1/2}.
\end{equation}
For all  $s>0$, for $m\geq 5$, one has
\begin{equation}\label{deux}
\left\Vert \frac{H_{t+s}F -H_sF}{t} -\frac{1}{2} 
\Delta H_sF\right\Vert_{m-5,\eps} \leq C_m ||F||_{m,\eps} t^{1/2}.
\end{equation}
\end{enumerate}
\end{prop}
{\it Proof.}
Formula (\ref{Tayloretendue}), integrated with respect to  $Y$ on $B^2$,
gives, for $k\leq m-1$:
\begin{equation}\label{Tayloretenduechaleur}
\begin{array}{lll}
\displaystyle 
\int_{B^2}
\widetilde{F}(X+Y) \ d\mu_{B^2,t}(Y)  = F(X) 
 \displaystyle+ \quad \sum_{i=1}^{k-1} \frac{1}{i!}
\sum_{J\in \Gamma^i, \delta \in \{0,1\}^i }
\int_{B^2}
\left( \prod_{r=1}^i \ell_{w_{j_r}^{\delta_r}}(Y)\right)
 \ d\mu_{B^2,t}(Y)
 \frac{\partial^i F}{\partial w_{j_1}^{\delta_1}\dots\partial w_{j_i}^{\delta_i}}(X)\\
  +\displaystyle
\sum_{J\in \Gamma^k, \delta \in \{0,1\}^k }
\int_{B^2}
\left( \prod_{r=1}^k \ell_{w_{j_r}^{\delta_r}}(Y) \right)
\int_0^1 \frac{(1-s)^{k-1}}{(k-1)!} 
 \ppp\left(\frac{\partial^k F}{\partial w_{j_1}^{\delta_1}\dots\partial w_{j_k}^{\delta_k}}\right)
(X+sY) \ ds     \ d\mu_{B^2,t}(Y) .
\end{array}
\end{equation}
We denote by $R_k$ the last term in the preceding formula. We have seen in 
subsection  \ref{subsection-Taylor-prolstoch}  that these functions 
admit $L^1$ norms, which allows us to exchange sums and integrals on $B^2$.
Using Wick's formula, we see that odd order terms are equal to $0$. One
can give a bound for the rest: 
\begin{lemm}
For all  $X\in H^2$, with $F\in S_m(\bbb,\eps)$ and $k\leq m-1$, one has
$$
|R_k(X)|\leq \frac{1}{\sqrt{\pi} \  k!}||F||_{m,\varepsilon}
 2^{\frac{3k}{2}} t^{\frac{k}{2}}
\Gamma(\frac{k+1}{2}) (\sum_{\Gamma} \varepsilon_j)^k.
$$
\end{lemm}
 {\it Proof.}
Notice that 
 $ \ppp\left(\frac{\partial^k F}{\partial w_{j_1}^{\delta_1}\dots\partial
 w_{j_k}^{\delta_k}}\right) $
  is bounded by
$||F||_{m,\varepsilon}\varepsilon_{j_1} \dots \varepsilon_{j_k} $.
One applies  H\"older's formula to the product of $\ell$ functions and one
 sums over $j_1,\dots, j_k$.\hfill $\square$

Even order terms allow us to find (thanks to Wicks formula) the successive
 powers of the Laplace operator and we get
\begin{equation}\label{machin}
 (H_tF)(X)=
\int_{B^2}
\widetilde{F}(X+Y) \ d\mu_{B^2,t}(Y)  = F(X)
 + \sum_{0<2p\leq k-1} \frac{1}{p!} \left(\frac{t}{2} \right)^p \Delta^pF(X) +R_k.
\end{equation}
Let us prove the point about  continuity. For $k=2$ and $m=3$, we can state
 the following
 result, since the rest is of order $t$:
$$
\forall X\in H^2,\quad
\left| H_tF(X) -F(X) -\right| \leq C_2 ||F||_{3,\eps} t,
$$
with $C_2= 2(\sum \eps_j)^{2}$.\\
This yields  (\ref{cont-0}) when $m=3$. To treat the general case one 
uses induction, working with $\partial^{\alpha}_u \partial^{\beta}_vF$, where
 $\alpha$ and $\beta$ have depth $1$ at most and using Proposition 
\ref{op-chaleur-H}. To obtain (\ref{contpartout}) one applies  $H_s$ to 
(\ref{cont-0}) (and loses one order of differentiability) and applies the
 semigroup property (Proposition \ref{semi-groupe!}).\\
Let us prove the point about differentiability. For $k=3$ and $m=4$, one can, 
in particular, obtain the following result since the rest is of order  $t^{3/2}$:
$$
\forall X\in H^2,\quad
\left| \frac{H_tF(X) -F(X)}{t} -
 \frac{1}{2} \Delta F(X)\right| \leq C_3 ||F||_{4,\eps} t^{1/2},
$$
with $C_3= \frac{1}{\sqrt{\pi} 3!} 2^{9/2} \Gamma(2) (\sum \eps_j)^{3}$.
This gives (\ref{un}) when $m=4$.  To treat the general case one 
uses induction, working with
 $\partial^{\alpha}_u \partial^{\beta}_vF$, where
 $\alpha$ and $\beta$ have depth $1$ at most and using Proposition 
\ref{op-chaleur-H}. To obtain (\ref{deux}) one applies  $H_s$ to (\ref{un}),
 (and loses one order of differentiability) and applies the
 semigroup property (Proposition \ref{semi-groupe!}).
\\ This completes the proof of Proposition \ref{auxiliaire}
\hfill $\square$

We now can state the main result about the heat operator in $S_m$ classes. For 
the sake of clarity, the two first points repeat former results of this
 subsection.
\begin{theo}\label{glavni-Sm}
Let $\varepsilon$ be  summable.
\begin{enumerate}
\item
For $m\geq 1$, the operator $H_t$  is continuous from  $S_m(\bbb,\varepsilon)$
to  $S_{m-1}(\bbb,\varepsilon)$ and for $m\geq 2$, the operator  $\Delta$ is
continuous from $S_m(\bbb,\varepsilon)$  to $S_{m-2}(\bbb,\varepsilon)$.
\item
For  $m\geq 3$, $H_t$ and $\Delta$ commute:
for all $F \in S_m(\bbb,\varepsilon)$,
$\Delta H_t F = H_t \Delta F \in S_{m-3}\bbb,(\varepsilon)$.
\item
Let $m\geq 6$ and $F\in S_m(\bbb,\eps)$. The application $t\mapsto H_tF$ is 
 $C^1$ from  $[0,+\infty[$ in $S_{m-6}(\bbb,\eps)$ and its derivative is
 $t\mapsto \frac{1}{2}H_t\Delta F$ .
\end{enumerate}
\end{theo}

{\it Proof.}
It remains to prove the last point. Set $\ph(t)= H_tF\in S_{m-1}(\bbb,\eps)$. 
According to the preceding proposition, $\ph$ is differentiable on
  $[0,+\infty[$ and $\ph'(t)= \frac{1}{2}\Delta H_tF=\frac{1}{2} H_t\Delta F$.
 But $H_t \Delta F\in S_{m-3}(\bbb,\eps) \subset S_{m-6}(\bbb,\eps) $. Since
 $\Delta F\in S_{m-2}(\bbb,\eps) $, an application of point $3$ (about 
continuity) proves that $t\mapsto H_t \Delta F$ is continuous from 
 $[0,+\infty[$ in $S_{m-6}(\bbb,\eps)$.
\hfill $\square$\\

\noindent{\bf Remark} It is not necessary to write $\Delta_{\bbb}$,
because of Remark  \ref{laplacien-independant}.
%%%%%%%%%%%%%%%%%%%%%%%%%%%%%%%%%%%%%%%%%%%%%%%%%%%%%%%%%%%%%%%
%%%%%%%%%%%%%%%%%%%%%%%%%%%%%%%%%%%%%%%%%%%%%%%%%%%%%%%%%%%%%%%
%%%%%%%%%%%%%%%%%%%%%%%%%%%%%%%%%%%%%%%%%%%%%%%%%%%%%%%%%%%%%%%
%%%%%%%%%%%%%%%%%%%%%%%%%%%%%%%%%%%%%%%%%%%%%%%%%%%%%%%%%%%%%%%
%%%%%%%%%%%%%%%%%%%%%%%%%%%%%%%%%%%%%%%%%%%%%%%%%%%%%%%%%%%%%%%
%%%%%%%%%%%%%%%%%%%%%%%%%%%%%%%%%%%%%%%%%%%%%%%%%%%%%%%%%%%%%%%
%%%%%%%%%%%%%%%%%%%%%%%%%%%%%%%%%%%%%%%%%%%%%%%%%%%%%%%%%%%%%%%
%%%%%%%%%%%%%%%%%%%%%%%%%%%%%%%%%%%%%%%%%%%%%%%%%%%%%%%%%%%%%%%

\subsection{The heat operator in the  classes $S(Q_A)$ } % 5.3
 In this subsection, the operator $A$ is self adjoint, nonnegative and trace
 class.
\begin{lemm}\label{diff-Ht-classeDD}
Let $f\in S(Q_A)$. For all  $m\in \N^*$ and all  $U_1,\dots, U_m$,
the application $ g_{m,U} : x\mapsto d^mf(x)(U_1,\dots, U_m) $ belongs to  $S(Q_A)$ and
$||g_{m,U}||_{Q_A} \leq ||f||_{Q_A} \prod_{j=1}^m Q(U_j)^{1/2}$.\\
The application $H_tf$ is differentiable on $H$ and 
$$
d(H_tf)(x) \cdot y = \int_{B} \ppp( u \mapsto df(u)\cdot y) (x+z) 
\ d\mu_{B,t}(z) =(H_tg_{1,y})(x).
$$
Moreover  
$$
|H_tf(x+y)-H_tf(x) -(H_tg_{1,y})(x)|\leq \frac{1}{2} ||f||_{Q_A} Q_A(y).
$$
\end{lemm}
{\it Proof.} 
One checks that $g_{m,U}$ is $C^{\infty}$ and that, for all integer $k\geq 1$ and
all  $h_1,\dots, h_k\in H$,
$$
d^kg(x)\cdot(h_1,\dots, h_k) = d^{m+k}f(x)\cdot (h_1,\dots, h_k, U_1,\dots, U_m).
$$
This proves that  $ g_{m,U} \in S(Q_A)$.\\
For $x,y\in H$, Taylor's formula gives 
$$
\tau_yf(x)= f(x) + df(x)\cdot y +\int_0^1 (1-s) d^2f(x+sy) \cdot y^2 \ ds.
$$
We denote by $R_2(x,y)$  the last term of this sum. Since  $\tau_y f, f$ and
 $x\mapsto df(x)\cdot y$ have stochastic extensions 
in  $L^p(B, \mu_{B,t})$, so does $x\mapsto R_2(x,y)$. One gets 
$$
H_tf(x+y) = H_tf(x)+ \int_B (\ppp( df(\cdot) \cdot y)(x+z)\ d\mu_{B,t}(z) 
+ \int_B \widetilde{R_2}(x+z) \ d\mu_{B,t}(z).
$$
One checks that the first integral gives a linear application with respect to
 $y$. The hypotheses on $f$ prove its continuity and the bound on the rest.
\hfill $\square$\\

By induction on the order $m$ on can deduce the following result:
\begin{prop}\label{5.10}
Let $f\in S(Q_A)$. For all $t>0$, the application $H_tf$ belongs to $S(Q_A)$
and  $||H_tf||_{Q_A} \leq ||f||_{Q_A}  $. Moreover, for all  integer $m$ and 
all $x,y_1,\dots, y_m$, one has, with the preceding notations,
$$
d^m(H_tf)(x)\cdot (y_1,\dots, y_m) = H_t(g_{m,y_1,\dots, y_m})(x).
$$
\end{prop}

We denote by $\Delta f(x)= {\rm Tr}(d^2f(x))$ the trace of the operator $M_x$ 
satisfying $<M_xU,V>= d^2f(x)(U,V)$ for all vectors $U,V$ of $H$. Its existence
 is ensured  by the inequalities  (\ref{Classe-DD-def-ineq}) and one can see
 it, too, as a sum of partial derivatives (with respect to an arbitrary
orthonormal basis of $H$). One can state the following proposition:
\begin{prop}\label{Htlaplaciencomm}
If $f\in S(Q_A)$, then $\Delta f \in  S(Q_A)$ with 
$ || \Delta f||_{Q_A}\leq {\rm Tr}(A) ||  f||_{Q_A}$. Moreover,  for all  $t>0$,
$$
\Delta(H_tf)(x)= H_t(\Delta f)(x).
$$
\end{prop}

{\it Proof.}
Let  $(e_j)$ be an orthonormal basis of  $H$. One can write
$$
{\rm Tr}(d^2f(x))= \lim_{n\rightarrow \infty}
\sum_{s=1}^n d^2 f(x)\cdot(e_s,e_s)
= \lim_{n\rightarrow \infty}  \sum_{s=1}^n g_{2,e_s,e_s}(x),
$$ 
with the notations of  Lemma \ref{diff-Ht-classeDD}.
Then the series $\sum g_{2,e_s,e_s}$ converges in $S(Q_A)$ because
$|| g_{2,e_s,e_s}||_{Q_A}\leq ||f||_{Q_A} <Ae_s,e_s>$ and $A $ is trace class.
Hence $\Delta f \in  S(Q_A)$ with
$ || \Delta f||_{Q_A}\leq {\rm Tr}(A) ||  f||_{Q_A}$.
Since  $H_t$ is continuous on $S(Q_A)$, one has
$$
{\rm Tr}(d^2H_tf(x))= \lim_{n\rightarrow \infty}
\sum_{s=1}^n d^2 H_tf(x)\cdot(e_s,e_s)
= \lim_{n\rightarrow \infty}  H_t(\sum_{s=1}^n g_{2,e_s,e_s})(x)= 
 H_t( \sum_{s=1}^{\infty} g_{2,e_s,e_s})(x)= H_t({\rm Tr}(d^2 f))(x).
$$
\hfill $\square$

\begin{prop}\label{chaleur-ClasseDD}
For all  $f\in S(Q_A)$, one has 
$$
\lim_{t\rightarrow 0}\left\Vert 
 \frac{H_t(f) -f}{t} -\frac{1}{2}\Delta f \right\Vert_{Q_A}=0.
$$
Moreover,  for all  $s>0$, one has
$$
\lim_{t\rightarrow 0} \frac{(H_{t+s}f) -H_sf}{t} =\frac{1}{2} 
 {\rm Tr}(d^2H_sf)=\frac{1}{2} \Delta H_sf=\frac{1}{2}  H_s\Delta f,
$$
the convergence taking place in  $S(Q_A)$.
\end{prop}

\hfill $\square$

{\it Proof.} 
 Let $x\in H$.  First prove that
\begin{equation}
\lim_{t\rightarrow 0} \frac{(H_t(f))(x) -f(x)}{t} =\frac{1}{2} {\rm Tr}(d^2f(x))
=\frac{1}{2}\Delta f(x).
\end{equation}
For $y\in H$, Taylor's formula gives 
$$
f(x+y)= f(x)+\sum_{j=1}^k \frac{1}{j!} d^jf(x) \cdot y^j 
+\int_{0}^1 \frac{(1-s)^k}{k!} d^{k+1}f(x+sy) \cdot y^{k+1} \ ds.
$$
We denote by  $R_k(y)$ the last term of the sum just above. According to
 Remark \ref{rem-ext-trans}, $\tau_x f$ has a stochastic extension 
$\tau_x \tilde{f}$ in $L^p(B,\mu_{B,h})$, with respect to the variable $y$.
 Indeed, $f$ admits a stochastic extension for all $p$ and Definition 
 \ref{Classe-DD-def}  implies that it is Lipschitz continuous. 
 By substraction, the rest
  $R_k$  also admits a stochastic extension $\widetilde{R_k}$. 
This extension  is bounded as follows: 
\begin{lemm}\label{majoration-Rtilda}
Let $t>0$, $p\in [1,+\infty[$. For $k\in \N^*$, one has
$$
||\widetilde{R_k}||_{L^p(B,\mu_{B,t})}\leq \frac{1}{(k+1)!} ||f||_{Q_A} C(p(k+1))^{k+1}
 S^{\frac{k+1}{\alpha(p(k+1))}} t^{\frac{k+1}{2}}, 
$$ 
with  $S=\sum_j \lambda_j$.
\end{lemm}
{\it Proof.} Let $(E_n)_n$ be an increasing sequence  of $\fff(H)$,
whose union is dense in $H$. Then
$$
\begin{array}{lll}
\displaystyle
||\widetilde{R_k}||_{L^p(B,\mu_{B,t})}
& \displaystyle\leq ||\widetilde{R_k}-R_k\circ \tilde{\pi}_{E_n} ||_{L^p(B,\mu_{B,t})}
+||R_k\circ \tilde{\pi}_{E_n} ||_{L^p(B,\mu_{B,t})}\\
& \displaystyle\leq ||\widetilde{R_k}-R_k\circ \tilde{\pi}_{E_n} ||_{L^p(B,\mu_{B,t})}
+||f||_{Q_A} \ || Q_A^{\frac{k+1}{2}}\circ \tilde{\pi}_{E_n}  ||_{L^p(B,\mu_{B,t})}\\
\end{array}
$$
by definition of $R_k$. Remark \ref{rem-maj-Q} enables us to give an upper
 bound independent of $n$ for the second term and to let $n$ 
converge to infinity. \hfill $\square$\\

One can then write, extending in $L^1(B,\mu_{B,t})$, according to 
 Proposition \ref{prol-diff-DD} :
 $$
\int_B \tilde{f}(x+y) d\mu_{B,t}(y)
 = f(x)+\sum_{j=1}^k \int_B \ppp\left( y\mapsto
\frac{1}{j!} d^jf(x) \cdot y^j \right)  d\mu_{B,t}(y) 
+ \int_B \widetilde{R_k}(y)  d\mu_{B,t}(y),
$$
where $\ppp$ represents the passage to the stochastic extension. For 
$j\leq k$ one uses the $L^1$ convergence and formula
(\ref{7AJN}) to obtain 
$$
\begin{array}{lll}
\displaystyle
 \int_B \ppp\left( y\mapsto
 d^jf(x) \cdot y^j \right)  d\mu_{B,t}(y)
 & \displaystyle= \lim_{n\rightarrow \infty} 
 \int_B 
 d^jf(x) \cdot \tilde{\pi}_{E_n}(y)^j d\mu_{B,t}(y) \\
&\displaystyle
 = \lim_{n\rightarrow \infty} 
 \int_{E_n} 
 d^jf(x) \cdot z^j d\mu_{E_n,t}(z) ,\\
\end{array}
$$
where $(E_n)_n$ is an increasing sequence  of $\fff(H)$, whose union is dense
 in $H$. For odd $j$, the terms are equal to $0$. For even $j$, one takes
an arbitrary orthonormal basis of $E_n$,  $(e_s)_{1\leq s\leq \dim(E_n)}$,
and one checks that 
$$
 \int_{E_n}  d^2f(x) \cdot z^2 d\mu_{E_n,t}(z)=
\sum_{s=1}^{\dim(E_n)} t \frac{\partial^2 f}{\partial e_s^2}(x).
$$
One then gets that, for any orthonormal basis of $H$, 
$$
\int_B \ppp\left( y\mapsto
 d^2f(x) \cdot y^2 \right)  d\mu_{B,t}(y) = t \sum_{j\in \N} 
 \frac{\partial^2 f}{\partial e_s^2}(x) = t {\rm Tr}(d^2f(x)).
$$
Applying the former reasoning to $k=3$ and using the upper bound of 
$\widetilde{R_3}$ in $L^1$ yield
$$
\left\vert 
\frac{(H_t(f)(x)-f(x))}{t}  -\frac{1}{2} {\rm Tr}(d^2f(x))
\right\vert \leq ||f||_{Q_A} \frac{1}{4!} C(4)^4 S^{\frac{4}{\alpha(4)}} t,
$$ 
which holds for all $x\in H$. This proves Formula (\ref{auxiliaire}).
Replacing $f$ by $g_{m, y_1,\dots ,y_m}$ in this inequality, we obtain, thanks to
 Lemma \ref{diff-Ht-classeDD} and Proposition \ref{5.10},
 $$
\begin{array}{ccc}
\displaystyle
\left\vert 
\frac{(d^mH_t(f)(x)\cdot(y_1,\dots, y_m) -d^mf(x)\cdot(y_1,\dots, y_m))}{t}
  -\frac{1}{2} d^m{\rm Tr}(d^2f(x))\cdot(y_1,\dots, y_m)
\right\vert\\
\displaystyle \leq
 || g_{m, y_1,\dots ,y_m} ||_{Q_A} \frac{1}{4!} C(4)^4 S^{\frac{4}{\alpha(4)}} t
\\ \displaystyle \leq
|| f ||_{Q_A} \prod Q_A(y_i)^{1/2}  \frac{1}{4!} C(4)^4 S^{\frac{4}{\alpha(4)}} t.
\end{array}
$$ 
One then has
$$
\left\Vert 
\frac{H_tf-f}{t}  -\frac{1}{2} \Delta f
\right\Vert_{Q_A} \leq \frac{1}{4!} C(4)^4 S^{\frac{4}{\alpha(4)}} t  ||f||_{Q_A},
$$
which gives the convergence in  $S(Q_A)$.\\
According to Proposition \ref{5.10}, $H_s$ is continuous on $S(Q_A)$ and its
 norm is smaller than $1$. The  semigroup property 
(Proposition \ref{semi-groupe!})  gives
$$
\left\Vert 
\frac{H_{t+s}f-H_sf}{t}  -\frac{1}{2} H_s \Delta f
\right\Vert_{Q_A} \leq \frac{1}{4!} C(4)^4 S^{\frac{4}{\alpha(4)}} t  ||f||_{Q_A},
$$
which achieves the demonstration of Proposition \ref{chaleur-ClasseDD},
since $H_s$ and $\Delta$ commute.
 \hfill $\square$

\begin{lemm}
Let $f\in S(Q_A)$, $x\in H$ and let $(e_n)$ be an arbitrary orthonormal basis
 of $H$. One denotes by $\frac{\partial}{\partial x_j}$ the differentiation
 in the direction of  $e_j$.
For all  integer $j$ one sets
$$
(\Delta^j)f(x)= \lim_{n\rightarrow \infty} \left(
\sum_{i=1}^n \frac{\partial^2}{\partial x^2_i}  \right)^j f(x).
$$
 One has for all  $h>0$, 
$$
H_tf(x)= f(x) + \sum_{j=1}^N \frac{1}{j!} \left(\frac{t}{2}\right)^j
\Delta^j f(x) + \int_B\widetilde{R}_{2N+1}(y)\ d \mu_{B,t}(y),
$$
with the upper bound of Lemma \ref{majoration-Rtilda}.
\end{lemm}
{\it Proof.}
One reasons as in the preceding demonstration but one considers $k=2N+1$ 
instead of stopping at $k=3$. For even $j$ one has
$$
\int_{E_n}  d^jf(x) \cdot z^j d\mu_{E_n,t}(z)= 
\int_{E_n} j! \sum_{\alpha\in \N^{\dim(E_n)}, |\alpha|= j} \frac{1}{\alpha!} 
\frac{\partial^j f}{\partial z^{\alpha} } z^{\alpha}  d\mu_{E_n,t}(z)
$$
and the terms where a coordinate of the multiindex $\alpha$ is odd are equal
to $0$. The computation of the other terms gives the result, thanks to
the equality
$$
\int_{\R} y^{2p} d\mu_{\R,1}(y) = \pi^{-1/2} 2^p \Gamma(p+1/2).
$$
\hfill $\square$

As a corollary of  Propositions \ref{Htlaplaciencomm} and
 \ref{chaleur-ClasseDD}, one can state the following commutation result, 
which will be used later on to prove a covariance result.

\begin{prop}\label{comp-phi}
Let $\ph $  be linear, continuous on $H$ and such that 
 $\ph^* \ph = \ph \ph^*= {\rm Id}_H$. Let $A$ be a linear application satisfying
the hypotheses of Definition
 \ref{Classe-DD-def}.
For all  $f\in S(Q_A)$, one can write
\begin{equation}
\forall t\geq 0, \ (H_tf)\circ \ph= H_t(f\circ \ph).
\end{equation}
\end{prop}
{\it Proof.}
One verifies that  $f\circ \ph $ (denoted by $f_{\ph}$) is in  
$S(Q_{\ph^* A \ph})$, with
$$
d^2f_{\ph}(x)\cdot (U,V)= d^2f(\ph(x)) \cdot(\ph(U),\ph(V))= 
<\ph^* M_{\ph(x)}(f) \ph U,V>
$$
and
$$
||f_{\ph}||_{Q_{\ph^* A \ph}} = ||f||_{Q_A}.
$$
(We still denote here by  $\Delta f(x)= {\rm Tr}(d^2f(x))$ the trace of the 
operator $M_x$ satisfying $<M_xU,V>= d^2f(x)(U,V)$ for all vectors $U,V$ in 
$H$.)
Moreover, the operator  $\ph^* M_{\ph(x)}(f) \ph$ is trace class and has
the same trace as $M_{\ph(x)}(f)$. Thus
$$
\Delta (f_{\ph}(x))= {\rm Tr}(d^2 f_{\ph}(x)) = {\rm Tr}(\ph^* M_{\ph(x)}(f) \ph)=
 {\rm Tr}( M_{\ph(x)}(f))=  {\rm Tr}(d^2 f(\ph(x))) =(\Delta f)(\ph(x)).
$$ 
Applying \ref{chaleur-ClasseDD} and the above remark to  $f\circ\ph$, one gets
 that 
$$
\lim_{t\rightarrow 0} \frac{H_t(f_{\ph}) -f_{\ph}}{t}  
= \frac{1}{2}\Delta (f_{\ph})= \frac{1}{2}(\Delta f)\circ \ph 
  \quad {\rm in}\quad  S(Q_{\ph^*A\ph}).
$$
Composing with $\ph^*$, one obtains that 
$$
\lim_{t\rightarrow 0}\left( \frac{H_t(f_{\ph}) -f_{\ph}}{t}  \right)\circ \ph^*
= \frac{1}{2}(\Delta f)= \lim_{t\rightarrow 0}\left( \frac{H_t(f -f}{t}  \right) 
 \quad {\rm in}\quad  S(Q_{A}).
$$
If one denotes by $T_t$ the operator defined on  $S(Q_A)$ by
$T_tf = H_t(f\circ \ph)\circ \ph^*$, one can verify that $(T_t)$ is a 
 semigroup on  $S(Q_A)$. Since both semigroups  $(T_t)$ and $(H_t)$ have the
 same infinitesimal generator $\frac{1}{2} \Delta$, which is continuous on 
 $S(Q_A)$ (Proposition \ref{Htlaplaciencomm}), they are uniformly continuous
 and equal (\cite{P}, Theorems 1.2 and 1.3, Chapter 1). This achieves the proof.
\hfill $\square$\\

We now can state the main result of this part. For the sake of clarity, the 
 first points repeat former results of the same part.

\begin{theo}\label{pcpal-DD}
Let  $A$  be a linear application on $H$ satisfying the hypotheses of
 Definition \ref{Classe-DD-def}. For all $f\in S(Q_A)$,  $\Delta f\in S(Q_A)$
 with $ ||\Delta f||_{Q_A}\leq {\rm Tr}(A) ||f||_{Q_A}$. Moreover, 
 for all  $t>0$,
$
\Delta(H_tf)(x)= H_t(\Delta f)(x).
$\\
 The function $t\mapsto H_tf$ is $C^{\infty}$ on  $[0,\infty[$
with values in  $ S(Q_A)$, with
$$
\frac{d^m}{dt^m}H_tf = \left(\frac{1}{2}\Delta\right)^m H_tf.
$$
 For all  $N\in \N^*$, one has
$$
H_tf= f + \sum_{k=1}^N \frac{t^k}{k!} \left(\frac{1}{2}\Delta \right)^kf 
+ t^{N+1} R_N(t), 
$$
where $R_N\in S(Q_A)$  is bounded independently of $t\in[0,1]$. 
\end{theo}

{\it Proof.} The first point comes from Proposition \ref{5.10}.
For the differentiability, according to Proposition \ref{chaleur-ClasseDD},
the result holds for $m=1$. But then, since $\Delta$ commutes with $H_t$ 
(Proposition \ref{Htlaplaciencomm}), one concludes by induction on $m$.\\
For the second point, one applies one of Taylor's formulae to $t\mapsto H_t f$,
which gives
$$
||R_N(t)||_{Q_A} \leq \frac{1}{(N+1)!} \sup_{s\in [0,t]}\left\Vert H_s(
\left(\frac{1}{2}\Delta \right)^{N+1}f )\right\Vert_{Q_A}
 \leq \frac{1}{(N+1)!} \left\Vert
\left(\frac{1}{2}\Delta \right)^{N+1}f \right\Vert_{Q_A}
$$
according to  Proposition \ref{5.10}. \hfill $\square$

%%%%%%%%%%%%%%%%%%%%%%%%%%%%%%%%%%%%%%%%%%%%%%%%%%%%%%%%%%%%%%%%%%%%%%
%%%%%%%%%%%%%%%%%%%%%%%%%%%%%%%%%%%%%%%%%%%%%%%%%%%%%%%%%%%%%%%%%%%%%%
%%%%%%%%%%%%%%%%%%%%%%%%%%%%%%%%%%%%%%%%%%%%%%%%%%%%%%%%%%%%%%%%%%%%%%
%%%%%%%%%%%%%%%%%%%%%%%%%%%%%%%%%%%%%%%%%%%%%%%%%%%%%%%%%%%%%%%%%%%%%%

\section{Appendix}\label{Annexe}% 7
We list here very general results used in the main part of the article.

\begin{lemm}\label{Hoelder-telescopique}
Let $(\Omega,\ttt,m)$ be a measure space. Let $N\geq 2$ be an integer.
For  $i\leq N$, let $f_{i}, g_i$  be functions on $\Omega$  with values in 
 $\R$ such that, for all $ p\in [1,+\infty[$, 
$ \ f_{i} \in L^p(\Omega,\ttt,m),\ g_i 
\in L^p(\Omega,\ttt,m).
$
For all  $p\in [1,+\infty[,$ set
$
M_p= \max_{1\leq i\leq N}( ||g_i||_{p} ,  || f_{i}||_{p}). 
$
Then  for all  $p \in [1,+\infty[,$ 
\begin{equation}\label{HT}
\left \Vert \prod_{i=1}^N f_{i} - \prod_{i=1}^N g_{i}  \right\Vert_p
\leq (M_{pN})^{N-1} \sum_{k=1}^N || f_{k}-g_k||_{pN}.
\end{equation}
More precisely, one has
\begin{equation}\label{HT-precise}
\left\Vert \prod_{i=1}^N f_{i} - \prod_{i=1}^N g_{i} \right\Vert_p
\leq \sum_{k=1}^N
   \left( \prod_{i=1}^{k-1}\left\Vert  f_{i}\right\Vert_{pN}
\prod_{i=k+1}^N\left\Vert  g_{i} \right\Vert_{pN}\right)
\left\Vert (f_{k}-g_k)\right\Vert_{pN}.
\end{equation}

\end{lemm}

{\it  Proof.}
According to the general  H\"older's inequality, the products belong to every 
$L^p$ since
$|| \prod_{i=1}^N g_i||_p \leq \prod_{i=1}^N || g_i||_{pN}$ (for example).
One decomposes $\prod_{i=1}^N f_{i} - \prod_{i=1}^N g_{i} $ as
$$
\begin{array}{lll}
\displaystyle
\prod_{i=1}^N f_{i} - \prod_{i=1}^N g_{i} &
\displaystyle
=\sum_{k=1}^N \left( \prod_{i=1}^k f_{i}\prod_{i=k+1}^N g_{i}
-\prod_{i=1}^{k-1} f_{i}\prod_{i=k}^N g_{i} \right)
& \displaystyle
=\sum_{k=1}^N \left( \prod_{i=1}^{k-1} f_{i}\prod_{i=k+1}^N g_{i} \right)
(f_{k}-g_k)\\
\end{array}
$$
agreeing that a product on an empty set of indices is equal to $1$.
 H\"older's inequality gives
$$
\left\Vert \left( \prod_{i=1}^{k-1} f_{i}\prod_{i=k+1}^N g_{i} \right)
(f_{k}-g_k)\right\Vert_p \leq 
  \prod_{i=1}^{k-1}\left\Vert  f_{i}\right\Vert_{pN}
\prod_{i=k+1}^N\left\Vert  g_{i} \right\Vert_{pN}
\left\Vert (f_{k}-g_k)\right\Vert_{pN}.
$$
The first $N-1$ factors are smaller than  $M_{pN}$. Taking the sum
gives the result.\hfill $\square$\\

\begin{coro}\label{corHT}
Let $F_1,\dots, F_N$ $N$ be functions defined on $H^2$ and admitting stochastic
 extensions $\widetilde{F_1},\dots, \widetilde{F_N}$ in  $L^p(B^2,\mu_{B^2,h})$ 
for all  $p\in[1,+\infty[$. Then $\prod_{i=1}^NF_i$ admits
$\prod_{i=1}^N\widetilde{F_i}$  as a stochastic extension  in
  $L^p(B^2,\mu_{B^2,h})$ for all  $p\in[1,+\infty[$.
\end{coro}
 {\it  Proof.}
Let  $(E_n)$ be  an increasing sequence   of $\fff(H^2)$, whose union is 
dense in $H^2$. According to (\ref{HT}),
$$
|| \prod_{i=1}^NF_i\circ \tilde{\pi}_{E_n} - 
\prod_{i=1}^N\widetilde{F_i}||_{p}
\leq \left( \sup_{n\in \N, i\leq N} (||F_i\circ \tilde{\pi}_{E_n} ||_{Np},
|| \widetilde{F_i} ||_{Np})\right)^{N-1} \sum_{i=1}^{N} 
|| F_i\circ \tilde{\pi}_{E_n}-\widetilde{F_i} ||_{Np},
$$
which gives the result.  \hfill $\square$

\end{document}